\documentclass[10 pt]{article}
\usepackage{amsmath,amsthm,amssymb,amscd}
\theoremstyle{plain}

\newtheorem{thm}{Theorem}[]
\newtheorem{pro}[thm]{Proposition}
\newtheorem{lem}[thm]{Lemma}
\newtheorem{cor}[thm]{Corollary}
\newtheorem{con}[thm]{Conjecture}

\theoremstyle{definition}

\newtheorem{dfn}[thm]{Definition}
\newtheorem{conv}[thm]{Convention}
\newtheorem{exa}[thm]{Example}
\theoremstyle{remark}

\newtheorem{rem}[thm]{Remark}

\DeclareMathOperator{\rk}{rk}

\DeclareMathOperator{\ad}{ad}
\DeclareMathOperator{\Hom}{Hom}

\DeclareMathOperator{\Br}{Br}
\DeclareMathOperator{\Bs}{Bs}

\DeclareMathOperator{\Pic}{Pic}
\DeclareMathOperator{\Gal}{Gal}
\DeclareMathOperator{\Ins}{Int}

\DeclareMathOperator{\NE}{\overline{NE}}
\DeclareMathOperator{\NA}{\overline{NA}}

\DeclareMathOperator{\NM}{\overline{NM}}

\renewcommand{\o}{\mathcal{O}}
\renewcommand{\l}{\mathcal{L}}
\newcommand{\q}{\mathbb{Q}}
\renewcommand{\r}{\mathbb{R}}
\renewcommand{\a}{\mathbb{A}}
\newcommand{\p}{\mathbb{P}}
\newcommand{\z}{\mathbb{Z}}
\newcommand{\F}{\mathbb{F}}
\renewcommand{\c}{\mathbb{C}}
\newcommand{\cstar}{\c^\times}

\begin{document}
\bibliographystyle{alpha}

\title{Birational Geometry of 3-fold Mori Fibre Spaces}

\author{
 Gavin Brown\thanks{gavinb@maths.warwick.ac.uk}\\
 Mathematics Institute\\
 University of Warwick\\
 Coventry, CV4 7AL, UK
\and
 Alessio Corti\thanks{a.corti@dpmms.cam.ac.uk}\\
 Dep.~of Pure Math.~and Math.~Stat.\\
 University of Cambridge\\ 
 Cambridge, CB3 0WB, UK\\
\and
 Francesco Zucconi\thanks{zucconi@dimi.uniud.it}\\
 DIMI, Universit\`a di Udine\\
 via delle Scienze 208\\
33100 Udine, Italy}

\date{15 July, 2003}

\maketitle 

\begin{abstract}
In this paper we study standard 3-fold conic bundles $X/\p^2$ and stable
cubic del Pezzo fibrations $X/\p^1$. These are the key examples of 
3-fold Mori fibre spaces (Mfs). 
We begin systematically to chart the \emph{geography} of these
varieties, namely the classification of their deformation families 
and discrete invariants. 
Given a 3-fold Mfs $X/S$, we aim to understand the set of all Mfs
$Y/T$ with $X$ birational to $Y$. For example, we say that $X/S$ is
\emph{birationally rigid} if this is a 1-element set.
We state conjectures on the rigidity of Mfs, with the goal of 
approaching optimal criteria.
Our main results are summarised in Tables~\ref{tab:1} and~\ref{tab:2}
and Figure~\ref{fig:geog}.
We want this to become your standard guidebook to Mori fibre spaces
and their birational geometry. 
\end{abstract}
\pagebreak

\tableofcontents

\section{Introduction}
\label{sec:introduction}

\subsection{Mori fibre spaces}

A \mbox{3-fold} Mori fibre space---we often use the acronym
Mfs---is a \mbox{3-fold} $X$ together with an extremal
contraction $f \colon X \to S$ of \emph{fibering type}.
This means that $X$ has $\q$-factorial terminal
singularities, $-K_X$ is ample on fibres, the relative Picard rank 
$\rho = \rk N^1(X) -\rk N^1 (S)$ equals $1$, and $\dim S < \dim X$. This
structure is the higher-dimensional generalisation of minimally ruled
surfaces and $\p^2$. 

This paper is devoted to \emph{strict} \mbox{3-fold} Mori fibre
spaces, that is, we always assume $\dim S\geq 1$. There are two cases:
\emph{conic bundles} when $\dim S =2$ and \emph{del
  Pezzo fibrations}---we write $dP_k$ fibrations when fibres are del
Pezzo surfaces of degree $k$---when $\dim S =1$. We also always assume that
the base variety $S$ is rational.

\begin{dfn} Let $f \colon X\to S$, $g \colon Y \to T$ be Mori fibre spaces. 
A birational map
\[
\varphi \colon X \dasharrow Y
\]
is \emph{square} if it maps a general fibre of $f$
isomorphically to a general fibre of $g$. We say that $X\to S$
and $Y \to T$ are \emph{square birational}, or \emph{square
  equivalent}, when there is a square birational map between them. 
\end{dfn}

\subsection{Geography}

Our first goal in this paper is to begin a systematic study of the
\emph{geography} of strict Mori fibre spaces. We focus on the case of
\emph{standard} conic bundles $X\to \p^2$ and \emph{stable} $dP_3$
fibrations $X\to \p^1$. 

Recall that a conic bundle $X \to S$ is standard if $X$ is
nonsingular; this in turn implies that $S$ is also nonsingular, and
the discriminant locus $\Delta \subset S$ is a nodal curve. Every
conic bundle is square birational to a standard conic bundle (it may
be necessary to blow up $S$). With qualifications, the assumption is not really
restrictive. In Section~\ref{sec:conic-bundles} of this paper we begin
a systematic study of standard conic bundles $X\to \p^2$. Necessarily,
$X$ embeds in a projectivised rank 3 vector bundle on $\p^2$, and we
learn how to write down the bundle and equations of $X$ explicitly. As
an illustration, we construct standard conic bundles for all
ramification data with a discriminant curve of degree 7. (There are 4
deformation families.) Our methods apply to higher degrees and to
surfaces other than $\p^2$; what we do here is a starting point
for future work.

The definition of stable del Pezzo fibration is technical; in this
paper we almost always work with the slightly stronger condition that
$X$ be nonsingular. Again, a $dP$ fibration is always square
equivalent to a stable one. 
In Section~\ref{sec:dp_3-fibrations} of this paper we
begin a systematic study of the geography of stable $dP_3$ fibrations 
$X\to \p^1$. Necessarily, $X$ is a member of a linear system $|3M+nL|$
on a rational scroll $\F(0,a,b,c)$. We determine the invariants
$n,a,b,c$ such that a general element of the linear system is a stable
$dP_3$ fibration. We study the degree 3 case because it is the most
interesting, but our methods apply to degree 1 and 2. 
Grinenko has studied del Pezzo fibrations of degree 1 and 2 
in a series of recent papers \cite{grinenko:01a, grinenko:01b, 
grinenko:01c, grinenko:00, gri02}, but he always works under the
assumption that the total space $X$ is Gorenstein. In the case of $dP$
fibrations of degree 1 and 2, this assumption is too restrictive. 

\subsection{Sarkisov links}
\label{sec:sarkisov-links}

Both in the case of conic bundles and del Pezzo fibrations $X/S$, we study
in detail nonsquare Sarkisov links originating from $X/S$. 
Recall that Sarkisov links are
the ``elementary'' building blocks of the birational geometry of Mfs:
any birational map between Mfs can be factored as a chain of Sarkisov
links. We usually, but not always, assume that $X$ is nonsingular and
sufficiently general in moduli.   
To understand Sarkisov links, we are led to compute explicitly the 
Mori cone of $X$ and the cone of mobile divisors on $X$. In our
examples $X$ has rank 2, so these cones 
are 2-dimensional, and we just need to determine
the two edges in terms of explicit loci on $X$. The mobile cone of $X$
is partitioned into chambers (in our case these are just wedges)
corresponding to moves of the 2-ray game \cite{corti:00} starting with $X$.
We study these moves explicitly to see if the 2-ray 
game leads to a Sarkisov link from $X\to S$ to a new Mori fibre space.
When $X$ is embedded in a scroll $\F$ over $S$ the moves are
often, but not always, induced by moves of $\F$.

\subsection{Rigid Mori fibre spaces and known criteria}

\begin{dfn} \cite[Foreword]{MR2001f:14004}
  The \emph{pliability} of a Mori fibre space $X \to S$ is the set
\[
 \mathcal{P}(X/S) = \{ \text{Mfs} \; Y/T \mid Y \; \text{is
 birational to} \quad X\}/\sim
\] 
 where $\sim$ denotes \emph{square birational equivalence}. We say
 that $X\to S$ is \emph{birationally rigid} when $\mathcal{P} (X/S)$
 consists of one element.
\end{dfn}

\noindent A few general criteria for birational rigidity have been known
for quite some time, see \cite{MR98e:14009}, \cite{corti:00}.
These criteria and their proofs are based on exploiting
properties of the 1-cycle $K^2_X$. The importance of the following
condition was first recognised in a series of brilliant papers by 
A.\ Pukhlikov. His point of view is explained for example 
in \cite{MR2001j:14010}.

\begin{dfn} \cite{MR2000d:14017, MR99f:14016, MR98j:14014} 
We say that a variety $X$ satisfies the \emph{$K^2$~condition} if
$K_X^2$ is not in the interior of the Mori cone $\NE X$ of effective
1-dimensional cycles on $X$.  
\end{dfn}

\begin{thm} \cite{MR84h:14047, MR82g:14035, MR81d:14023},
\cite{MR2000d:14017, MR99f:14016, MR98j:14014} 
Let $X$ be a nonsingular \mbox{3-fold} and $X\to S$ a conic bundle or a $dP_3$
fibration (satisfying an additional ``genericity'' condition on the
singularities---see Theorem~\ref{thm:dp3} below). If 
$X$ satisfies the $K^2$~condition, then $X \to S$ is birationally
rigid.  
\end{thm}

More recent results explore examples that lie on the boundary between rigid and
nonrigid \cite{cm02, mella02}, \cite{grin03}, \cite{MR2003g:14017}.
These papers study Mfs $X/S$  which are either \emph{bi-rigid} in the 
sense that they have exactly two models as Mori fibre spaces 
(more precisely, $\mathcal{P}(X/S)$ is a set with two elements), or that are
rigid in a subtle way, that is, there exist Mfs $Y/T$ and nonsquare maps $X
\dasharrow Y$, but $Y/T$ is always square birational to $X/T$.
In this paper we construct many examples of varieties that are likely
to behave similarly. 

\subsection{Conjectures}
\label{sec:mfs_con}

The known general criteria for rigidity of Mori fibre spaces are not
optimal. Our second goal in this paper is to state some conjectures, based on
the intuition built on looking at the examples we construct. Some of
these conjectures we expect to be able to prove in the near future,
others, we feel, are definitely harder, and some are probably wrong. 
The $K^2$~condition is not natural. Our conjectures are stated in term
of a geometrically more natural condition, whose importance was also
recognised in work by Grinenko \cite[Conjecture 1.5]{grin99}
\cite[Conjecture 1.6]{grinenko:00} \cite[Conjecture 2.5]{gri02}.

\begin{dfn} \label{dfn:con-star} A variety $X$ \emph{satisfies 
condition~$(\ast)$} if the anticanonical class $-K$ is not in the interior of
  the cone of mobile divisors of $X$.
\end{dfn}

\noindent Note that condition~$(\ast)$ is stronger than the $K^2$~condition.
Indeed, if $-K$ is in the interior of the mobile cone, we can write
for some $n>0$
\[
-nK=M_1+M_2=H_1+H_2
\]
where $M_1\sim H_1$ have no component in common, and,
similarly, $M_2\sim H_2$ have no component in common.
Then writing
\[
n^2K^2=\sum M_iH_j
\]
shows that $K^2$ is in the interior of the Mori cone.

In Sections~\ref{sec:rigid-conic-bundles} and~\ref{sec:dp-con}
we state some conjectures on the rigidity of standard conic bundles
and stable $dP_3$ fibrations in terms of condition~$(\ast)$. 
For many of the examples that we construct which seem to lie on the
boundary between rigid and nonrigid, it would be natural to conjecture
that they are bi-rigid, or that they are still rigid, though they do not
satisfy condition~$(\ast)$. We have resisted the temptation to make
such conjectures, but we invite our readers to make their own and
possibly prove some theorems. 

It is by now well understood how to use Pukhlikov's $K^2$~condition 
to prove rigidity of Mfs, but it seems difficult to go further.
On the other hand, we have as yet no experience or success
using condition~$(\ast)$ to prove rigidity of Mfs.  
We hope that this paper will help to introduce a 
new point of view on the birational geometry of Mori fibre spaces,
where the geometric properties of $-K$ play more prominent a role through
condition~$(\ast)$. 

\subsection{Tables and Figures}
\label{sec:tables-figures}

Tables~\ref{tab:1} and~\ref{tab:2}, and Figure~\ref{fig:geog}
summarise and collect key information which is obtained in many
calculations throughout the paper. 
Table~\ref{tab:1} shows the list of standard conic bundles over $\p^2$ 
with discriminant of degree 7 and contains information about their
nonsquare Sarkisov links to alternative models as Mori fibre spaces. 
Figure~\ref{fig:geog} depicts the geography of $dP_3$ fibrations and
contains information about the general member of each
family, showing whether it satisfies condition~$(\ast)$ and whether we know an
alternative model as a Mori fibre space.
This information is supplemented in Table~\ref{tab:2}, which lists all
families of $dP_3$ fibration for which we know an element that has a
nonsquare Sarkisov link to some Mfs $Y/T$ and contains detailed
information about the link. 

\subsection{The Appendix}
\label{sec:appendix}

The Appendix sets out our notation for rational scrolls and
information about their birational maps. A scroll $\F$ is a quotient of
affine space by an action of the group $G=\cstar \times \cstar$. Though
one is not normally aware of this, the quotient depends on the
choice of a $G$-linearisation. Different $G$-linearisations produce
different quotients. In this manner we can understand the 2-ray game
originating from the scroll. In the Appendix, we introduce notation for
this which is used throughout the paper when constructing Sarkisov
links of Mfs $X\subset \F$; indeed, it is often the case that
moves of $\F$ induce moves of $X$.

\subsection{Various cones}
\label{sec:various-cones}

Let $X$ be a projective variety. We denote $N^1(X)=N^1(X, \r)$ the
vector space of Cartier divisors on $X$ with real coefficients, modulo
\emph{numerical} equivalence, and $N_1(X)=N_1(X, \r)$ the dual space of
1-dimensional cycles on $X$ with real coefficients, modulo
numerical equivalence. Throughout this paper, it is crucial
to be aware of various cones associated to $X$:
\begin{description}
\item[The Mori cone] $\NE_1(X)$, or simply $\NE (X)$, is the closure
  of the cone in $N_1(X)$ generated by the effective 1-cycles.
\item[The ample cone] $\NA^1(X)$ is the cone in $N^1(X)$ generated by
  ample divisors; Kleiman's criterion states that it is the dual cone
  to the Mori cone.  
\item[The mobile cone] $\NM^1(X)$ is the cone in $N^1(X)$ generated by
  mobile divisors. We say that a divisor $D$ on $X$ is \emph{mobile}
  if a positive multiple $nD$ moves in a linear system with no fixed
  divisor.  
\item[The quasieffective cone] $\NE^1 (X)$ is the closure of the 
  cone in $N^1(X)$ generated by effective divisors. Divisor in this
  cone are called \emph{quasieffective}.
\end{description}

\subsection{Disclaimer}

This paper was prepared against a tight (for us) deadline. We did not
test some of the conjectures as much as we would have liked,
and it is likely that some of the definitions, calculations,
etc.\ contain mistakes. The responsibility for these is of course
ours, but we still want to apologise for the inconvenience we
may be causing our readers. 

\section{The Sarkisov category}
\label{sec:sarkisov-category}

The \emph{Sarkisov category} is the category of Mori fibre spaces and
birational maps between them. We begin with a simple result
stating that Mfs that have a model as a Fano \mbox{3-fold} belong to 
finitely many deformation families. 

\begin{pro} \label{pro:finite}
  There is a finite collection of algebraic families of Mori fibre
  spaces with the following property: if $X\to S$ is birational to a
  Fano \mbox{3-fold} $Y$, then $X/S$ is square birational to a member of one
  of the families in the collection.  
\end{pro}

\begin{proof}
  We briefly explain the idea of the proof. By the Sarkisov program
  \cite{MR96c:14013}, nonrigid Mfs are parametrised by nonsquare
  links of the Sarkisov program. These are in turn parametrised by
  weak terminal Fano \mbox{3-folds} of rank 2, and these are bounded by 
  \cite{MR2001h:14053}. 
\end{proof}

We hope that we will soon be able to prove the following. 

\begin{con}
  If $X\to S$ is a Mori fibre space, the pliability $\mathcal{P}
  (X/S)$ is in a natural way an algebraic variety. 
\end{con}

\noindent The key difficulty is this. Let $\mathcal{X}\to T$ be a family
of Fano \mbox{3-folds} with a birational selfmap 
$\sigma \colon \mathcal{X} \dasharrow \mathcal{X}$ that maps fibres
birationally into fibres and hence induces a birational selfmap of the
base $T$. If this induced birational selfmap of $T$ has
infinite order and $t_1\in T$ is a general point, then the fibre $X_{t_1}$  
is birational to $X_{t_2}$, and then to $X_{t_3}$ and
so on. It seems possible that all the values $t_i\in T$, but not the
whole of $T$, may contribute to the pliability.

Recent experience \cite{cpr} \cite{cm02, mella02} suggests the
following rather optimistic conjectures.

\begin{con}
  For the total space $X$ of a \mbox{3-fold} Mfs $X\to S$ to be rational is a
  topological property. In other words, if $X/S$ and $Y/T$ are \mbox{3-fold}
  Mori fibre spaces, $X$ is rational and $Y$ is diffeomorphic to $X$,
  then $Y$ is also rational.
\end{con}

\begin{con}
  For a \mbox{3-fold} Mori fibre space to be rigid is a topological
  property.  
\end{con}
\noindent In particular, being rigid or rational is constant along algebraic
families where all fibres are diffeomorphic. 

\begin{con}
  The pliability $\mathcal{P}(X/S)$ of a \mbox{3-fold} Mori fibre space is in
  a natural way a topological invariant.
\end{con}

\section{Conic bundles}
\label{sec:conic-bundles}

We do not remind you of the known sufficient conditions for rigidity of conic
bundles or Iskovskikh's conjectural characterisation of conic bundles
with rational total space 
\cite{MR97j:14044, MR93i:14033, MR92k:14036, MR88i:14038}
\cite{corti:00}. We do illustrate a
method to write down equations of a conic bundle with assigned
discriminant. This method is based on Catanese's work on the Babbage
conjecture \cite{MR83c:14026} and is explicitly computable, unlike the 
traditional
abstract approach via Brauer groups. We focus on conic bundles over
$\p^2$ with discriminant a nodal plane curve of degree~7. Motivated by
these examples, we state some conjectures on the rigidity of conic bundles.

\subsection{Conic bundles and Brauer groups}
\label{sec:conic-brauer}

For details on this section see \cite{MR48:299}, \cite{MR84h:14047}.
Let $X \to S$ be a \emph{standard} conic bundle over a rational surface
$S$. By definition, this means that $X$ is nonsingular
and the relative Picard rank is 1, a consequence of which is that
$S$ is itself nonsingular and the discriminant curve $\Delta \subset
S$ is nodal. Over
$\Delta$, the fibre is generically the sum of two lines, which
specifies a \mbox{2-to-1} \emph{admissible cover} $N \to \Delta$.
(When $C$ is nonsingular, admissible just
means \'etale; in general, an admissible \mbox{2-to-1} cover is required to
ramify over both branches of each node of $C$.)
Together, we refer to the double cover $N\to \Delta$ and the embedding $\Delta
\subset S$ as the \emph{ramification data} of the conic bundle.

It is well known that there is a standard conic bundle $X\to S$ with
any preassigned ramification data $\Delta \subset S$ and \mbox{2-to-1}
admissible cover $N \to \Delta$. The total space $X$ is not unique, but any two
choices are square birational over the base.

The traditional proof uses the exact sequence \cite{MR48:299} (in which $S$ is
a rational surface with field of rational functions $K$)
\[
0\to \Br S \to \Br K\to\bigoplus_{\text{curves}\;C \subset S} 
H^1(K(C), \q/\z) \to \bigoplus_{\text{points}\; P \in S} \mu^{-1} \to
\mu^{-1} \to 0
\]
and the fact, due to Platonov and reproduced in \cite{MR84h:14047},
that every element of order 2 in the Brauer group $\Br K$ can be
represented by a quaternion algebra. See the references for details
and explanations.

In the following subsections, we illustrate an effective proof of this 
statement in two steps: first we learn how to specify effectively a
\mbox{2-to-1} cover $N \to \Delta$, then we write down explicit equations for
a conic bundle $X$.

\subsection{Catanese's theory}
\label{sec:cataneses-theory}

Let $C$ be a nonsingular curve of degree $d$ in $\p^2$. (With a bit of
care, the theory works unchanged for nodal curves.) We fix coordinates
$u_0, u_1, u_2$ on $\p^2$ and denote $S=k[u_0,u_1,u_2]$ the
homogeneous coordinate ring. 
We are interested in admissible \mbox{2-to-1} covers of $C$. We summarise the
part of Catanese's paper \cite{MR83c:14026} which is relevant to us;
see also \cite{Dixon}.
The theory works identically for covers which in addition ramify over
a specified hyperplane section of $C$, and we treat these as well. 

An admissible \mbox{2-to-1} cover corresponds to a line bundle 
$\mathcal{L}$ on $C$ with a choice of isomorphism: either 
$\mathcal{L}^2 =\o_C$, or $\mathcal{L}^2 =\o_C(-1)$.
Denote by $i\colon C \hookrightarrow \p^2$ the inclusion. By 
\cite{MR83c:14026}, Theorems~2.16 and 2.19, $i_\ast \mathcal{L}$ is a 
Cohen-Macaulay $\o_{\p^2}$-module, so it has a 2-step symmetric locally 
free resolution:
\[
0 \to \oplus \o_{\p^2} (-l_i) \stackrel{A}{\to} \oplus \o_{\p^2} (-r_i) 
\to i_\ast \mathcal{L} \to 0.
\]
In other words, $A$ is a symmetric $n\times n$ matrix whose
entries are homogeneous forms on $\p^2$. 

Note that $\det A$ is a homogeneous equation of $C$. 
On the other hand, writing $C$ as a symmetric determinantal specifies the line
bundle $\l$ and the double cover. 

Denote by $d_i$ the degree of the $i$-th diagonal entry $a_{ii}$ of
$A$; $d = \sum d_i$ is a partition
of $d$. The degree of $a_{ij}$ is $(d_i+d_j)/2$,
so that all the $d_i$s have the same parity. 
We have the following numerology:
\begin{equation*}
\begin{aligned}
r_i &= (d+e-d_i)/2\\
l_j &= (d+e+d_j)/2
\end{aligned}
\quad
\text{where}
\quad
 e=
 \begin{cases}
   0 \quad \text{if}\quad  \l^2 = \o_C\\
   1 \quad \text{if}\quad  \l^2 = \o_C(1).
 \end{cases}
\end{equation*}
Let
\[
M = \bigoplus_{n \geq 0} M_n \quad \text{where} \quad M_n =
H^0 \bigl(\p^2, i_\ast \l (n)\bigr)
\]
be the Serre module of the sheaf $i_\ast \l$. 
The direct summand $\o_{\p^2} (-r_i)$ in the presentation
corresponds to a generator $m_i$ of degree $r_i$ of $M$.
The matrix $A$ encodes the structure of $M$ as a graded module over
the homogeneous coordinate ring $S=k[u_0,u_1,u_2]$ of $\p^2$, as well
as the identity $m_i m_j = b_{ij} \in S$  
in which $(b_{ij})=B = \ad A$ is the adjugate matrix of $A$.

We work out the well-known example
of theta characteristics of plane quartics in detail; this
corresponds to the case $d=4$, $e=1$ of the theory.
Here $\l$ is a line bundle with $\l^2 =\o(-1)$, and $\l(1)$ is a theta
characteristic since $K_C = \o(1)$. A plane quartic has genus $g=3$
and $\deg \l = -2$, so Riemann--Roch gives
\[
h^0(C,\l(1))-h^1(C,\l(1))=0.
\]
We have two cases: $ h^0(\l (1) )= 0$ (even theta characteristic), and
$ h^0(\l (1) )= 1$ (odd theta characteristic).  
We work out the odd case. 

By assumption, we have exactly one generator 
$m \in M_1 = H^0(C,\l(1))$ of degree 1. 
Next, $\l(2)$ has degree 6 and Riemann--Roch calculates 
\[
h^0(C,\l(2))=6+1-3=4,
\]
so $M_2$ is based by 4 elements: they are
$u_0m, u_1m, u_2m$
and a new generator $n$. Noting $h^0(C,\l(3)) = 8$, we see that there
must be a relation in degree 2 between the 9 elements $u_iu_jm$,
$u_kn$; this relation must take the form $n f_1 = m f_2$
where $f_1$ is a linear form in $u_0, u_1, u_2$ and $f_2$ a 
quadratic form. The matrix $A$ has the form
\[
A= 
\begin{pmatrix}
f_1& f_2 \\
f_2& f_3
\end{pmatrix}.
\]
The whole structure is specified by writing the curve
$C=(f_1f_3-f_2^2=0)$ as a symmetric determinantal.
This presentation reveals a preferred 
bitangent $f_1=0$, which ``is" the odd section $m$.

In the even case, you will find 4 generators $m_1,\dots,m_4$
of degree 2, and 4 linear relations
among the 12 elements $u_i m_j$: $A$ is a $4\times 4$
symmetric matrix of linear forms. 

\subsection{Plane curves of degree 7}
\label{sec:deg-7}
We study the case $d=7$, $e=0$, that is, line bundles $\l$ 
with $\l^2 = \o_C$ on plane curves $C$ of degree $\deg C = 7$.
Here $K_C = \o(4)$, so $\l(2)$ is a theta characteristic. A plane curve of
degree $7$ has genus $g=15$ and canonical degree $\deg K_C=28$. 
Riemann--Roch and Serre duality give at once
\begin{eqnarray*}
\label{eq:rr7}
h^0\bigl(C,\l(n)\bigr) - h^1\bigl(C,\l(n)\bigr) & = & 7(n-2)        \\
h^1\bigl(C,\l(n)\bigr) & = & h^0\bigl(C,K_C\otimes\l(n)^\ast\bigr)     \\
  & = & h^0\bigl(C,\l(4-n)\bigr).
\end{eqnarray*}
We use these equations without comment throughout.
It is clear that $p_1 = h^0(C,\l(1))$
and $p_2 = h^0(C,\l(2))$ determine all the $p_n=h^0(C, \l(n))$:  indeed,
the Hilbert series of $i_\ast\l$ as an $S$-module is
\begin{multline*}
P(t) =\sum_{n\geq 0} h^0(C, \l(n)) t^n = \\
 p_1t + p_2t^2 + (7 + p_1)t^3 + 14t^4 + \cdots + 7it^{i-2} + \cdots.
\end{multline*}

There are four cases, corresponding to the four possible values of
$p_2 = h^0(C, \l(2))$. We show later that the first two cases
occur on a generic curve, while the last two cases can only occur on
curves with special moduli.  In the four cases,
the degrees of the diagonal entries of the matrix $A$ are the
summands of the four different partitions of the number $d=7$
into odd summands.
\begin{enumerate}
\item $h^0(C, \l(2))=0$ (generic even theta),
\item $h^0(C, \l(2))=1$ (generic odd theta),
\item $h^0(C, \l(2))=2$ (special even theta),
\item $h^0(C, \l(2))=3$ (special odd theta).
\end{enumerate}

\subsubsection{$h^0(C, \l(2))=0$}
\label{sec:case-1}
We have $h^0(C, \l(3)) = 7$ generators of degree 3. Multiplying these
by $u_0$, $u_1$, $u_2$, we get 21 elements in the 14-dimensional space
$M_4$, so there are 7 equations between them.
In this case $A$ is a 7 by 7 symmetric matrix of linear forms.

\subsubsection{$h^0(C, \l(2))=1$}
\label{sec:case-2}
Again $h^0(C, \l(3))=7$, but here we assume a generator $m$ of
degree 2, giving $u_i m$ in $M_3$, so we need 4 more generators $n_1,\dots,n_4$
in $M_3$. In this case $A$ is a $5 \times 5$ symmetric matrix of homogeneous
forms of degrees
\[
\begin{pmatrix}
  3 & 2 & 2 & 2 &2 \cr
    & 1 & 1 & 1 &1 \cr
    &   & 1 & 1 &1 \cr
    &   &   & 1 &1 \cr
    &   &   &   &1
\end{pmatrix}
\]
in which the degrees of the diagonal elements express $d=7=3+1+1+1+1$
as a partition into odd summands.

\subsubsection{$h^0(C, \l(2))=2$}
\label{sec:case-3}

This case corresponds to the partition $7=3+3+1$. The assumption is
that there are two generators $m_1, m_2$ in $M_2$, so, as usual,
$h^2(C, \l(3))=7$ then shows that we need one more generator $n$
in degree 3. We conclude that $A$ is a $3 \times 3$ symmetric matrix
of homogeneous forms of degrees
\begin{equation}
\label{mx:case3}
\begin{pmatrix}
    3& 3& 2\cr
     & 3& 2\cr 
     &  & 1
\end{pmatrix}.
\end{equation}
We show below that not all plane septics have a theta-characteristic
of this type.  More precisely, the locus of those who do is of
codimension 1 in moduli.

\begin{exa}[Counting moduli]
Let $C = (\det A = 0)\subset\p^2$ for $A$ a symmetric
$3\times 3$ matrix of forms of degrees given by matrix
format (\ref{mx:case3}) above. The matrix A depends on 
$3\times 10 + 2\times 6 + 3 = 45$ parameters, and $A$ is uniquely 
determined by the choice of the generators of the
Serre module.  In this case, that choice amounts to the choice of
$m_1,m_2$ in $M_2$ (which has the 4 degrees of freedom of the basis of a 
2-dimensional vector space) and the choice of the additional generator
$n$ in $M_3$ (7 degrees of freedom) for a total of $4+7=11$ degrees of 
freedom. The space of septics of the given form has dimension $45 - 11 = 34$,
whereas the space of all septics has dimension $35$.  
\end{exa}

\subsubsection{$h^0(C, \l(2))=3$}
\label{sec:case-4}

This case corresponds to the partition $7=5+1+1$. Here 
$p_1=1$ and the element $m$ in $M_1$ generates a basis
$u_i m$ of $M_2$. This time $h^0(C,\l(3)) = 8$,
so we need two additional generators $n_1,n_2$ in degree 3.
The matrix $A$ is a $3\times 3$ matrix of homogeneous forms of degrees
\[
\begin{pmatrix}
    5& 3& 3\cr
     & 1& 1\cr
     &  & 1
\end{pmatrix}.
\]
Once again, the locus of septics who own a theta characteristic of
this type is of codimension 1 in moduli.

\subsection{Conic bundles over $\p^2$ with discriminant of degree 7}
\label{sec:delta-7}

We construct standard conic bundles over $\p^2$ with ramification data
a \mbox{2-to-1} admissible cover $N \to \Delta \subset \p^2$ of a plane curve
of degree $7$. We make four deformation families of such \mbox{3-folds}, 
following the analysis in the previous section. 
We show that general members $X$ of the first three families have an 
alternative model as a Mori fibre space. (This may be true for all
their members, but we didn't check that). Conjecture~\ref{con:cb} below
implies that a general member of the fourth family is birationally rigid.
This is a striking prediction, because this example is far
from the numerical range where the known rigidity criteria apply. 

\subsubsection{$h^0(\Delta, \l(2))=0$}
\label{sec:cb7-1}

Let $N \to \Delta$ be as in \ref{sec:case-1}. A conic bundle with this
ramification data was known classically, and it is birational to the
Fano \mbox{3-fold} $Y = Y_{2,2,2} \subset \p^6$, the codimension 3 complete
intersection of three quadrics in $\p^6$. Starting from $Y$, the
link $Y\dasharrow X$ begins with the projection from a line $\ell
\subset Y$. The details of this construction are well known, but see
also the following example.

\subsubsection{$h^0(C, \l(2))=1$}
\label{sec:cb7-2}
  Let $N \to \Delta$ be as in \ref{sec:case-2}.
  We want explicit equations for a conic
  bundle with ramification data $N \to \Delta$. 
  Consider the variety
 \[
  Z = \{ {}^t \mathbf{x} A \mathbf{x} = 0 \} \subset \F = \F (1,0,0,0,0). 
 \] 
 This is a bundle of 3-dimensional quadrics immersed in a scroll over
 $\p^2$ with ramification data $N\to \Delta$. Indeed, a singular
 3-dimensional quadric is, generically, a cone over $\p^1 \times
 \p^1$, and the two rulings give back the cover $N \to \Delta$. 
 It is natural to try to construct our conic bundle from $Z$ by a kind
 of ``dimensional reduction''. 

\paragraph{Claim 1}
\label{sec:claim-1}

 We can choose coordinates in the scroll $\F$ such that the matrix $A$
 has the form
\[
A=
\begin{pmatrix}
  B & \mathbf{b} \cr
  {}^t \mathbf{b} & 0 
\end{pmatrix}
\]
where $\mathbf{b}$ is a $1\times 4$ column vector of linear
forms. To see that this is possible, we look for a coordinate change
of the form:
\[
\begin{pmatrix}
  x_0 \\
\mathbf{\xi}
\end{pmatrix}
\mapsto
\begin{pmatrix}
  1            & 0 \\
\mathbf{\beta} & M
\end{pmatrix}
\begin{pmatrix}
  x_0 \\
\mathbf{\xi}
\end{pmatrix}
\]
where ${}^t\mathbf{\xi} = (x_1,x_2,x_3,x_4)$, $M=(m_{ij}) \in
GL(4,\c)$, and $\mathbf{\beta}$ is a $4\times 1$ column vector of
linear forms in $u_0,u_1,u_2$. The coordinate change
 brings $A$ in the 
wanted form if and only if 
the column vector with entries
$(0,m_{14},m_{24},m_{34},m_{44})$
is a nonzero null-vector of $A$. It is easy to see that such a vector
exists; indeed the $4\times 4$ symmetric submatrix $(a_{ij})$
of $A$ is a matrix of linear forms and corresponds to a net of
3-dimensional quadrics. A base point of the net of quadrics gives the
null-vector we want.

\paragraph{Equations of $X$}
\label{sec:equations-x}

Define
\[
X=\{ {}^t \mathbf{x} B \mathbf{x} = \mathbf{b} \cdot \mathbf{x}= 0 \}
\subset \F (1,0,0,0)
\]
where $\F(1,0,0,0)$ is the subscroll $\{x_4=0\} \subset \F(1,0,0,0,0)$.
One can easily check that $X$ is a conic bundle with the required 
ramification data $N \to \Delta$. Notice that, in terms of
projectivised vector bundles on $\p^2$, the linear form expresses $X$
as a conic in a projectivised nonsplit rank 3 vector bundle on $\p^2$. 

\paragraph{Claim 2}
\label{claim:claim-2}

If $\mathbf{b}$ is sufficiently general, we may further change
coordinates so that
\[
\mathbf{b}=
\begin{pmatrix}
  0 \\
  u_0\\
  u_1\\
  u_2
\end{pmatrix}.
\]

\paragraph{Model as a Fano \mbox{3-fold}}
\label{sec:model-as-fano}

We construct a Sarkisov link from $X$ to a 
codimension 3 Fano \mbox{3-fold} $Y=Y_{2,3,3,3,3} \subset
\p(1^6,2)$. The ideal of $Y$ is generated by
the $4\times 4$ Pfaffians of a $5\times 5$ antisymmetric matrix of homogeneous
forms of degrees
\[
\begin{pmatrix}
  2&2&2&2\\
   &1&1&1\\
   & &1&1\\
   & & &1
\end{pmatrix}
\]
on $\p(1^6,2)$. Assuming that $X$ is generic, we exhibit a birational
map, in fact a link of the Sarkisov program, from $X$ to a Fano \mbox{3-fold}
of the type just described. 
This alternative model of $X$ as a Mori fibre space
demonstrates, in particular, that $X$ is not birationally rigid.

If we choose coordinates as in Claim~\ref{claim:claim-2}, then $X \subset
\F(1,0,0,0)$ is given by equations
\[
\begin{cases}
  {}^t \mathbf{x} B \mathbf{x} =0 \\
  u_0x_1+u_1x_2+u_2x_3 = 0
\end{cases}
\]
where $B$ is a $4\times 4$ symmetric matrix of forms of degrees
\[
\begin{pmatrix}
  3 & 2 & 2& 2\\
    & 1 & 1& 1\\
    &   & 1& 1\\
    &   &  & 1
\end{pmatrix}
\]
in the coordinates $u_0,u_1,u_2$. 

The morphism $p\colon \F \to \p^5=\p H^0(\F, M)^\ast$ given by
\[
(x_1,x_2,x_3, x_4, x_5,x_6) =
(x_1,x_2,x_3,u_0x_0,u_1x_0,u_2x_0)
\]
identifies $\F$ with the blow up of $\p^5$ along the 2-plane $\Pi
=\{ x_4=x_5=x_6=0\}$. Under this morphism $X$ maps to
$p(X)=\overline{Y}_{2,3}$, a generic complete intersection of a
quadric and a cubic in $\p^5$ containing the plane $\Pi$.
To find the equations of $\overline{Y}$, note that we 
can find unique homogeneous quadrics
$q_0$, $q_1$, $q_2$ in six variables, such that
\[ 
{}^t \mathbf{x} B \mathbf{x} =\sum_{i=0}^2
u_iq_i(x_1,x_2,x_3,u_0x_0,u_1x_0,u_2x_0).
\]
(If $B$ is generic, the quadrics $q_0,q_1,q_2$ are also
generic.) It is easy to see that $\overline{Y}$ is given by the
following equations:
\[
\begin{pmatrix}
  x_1 & x_2 & x_3\\
  q_0 & q_1 & q_2
\end{pmatrix}
\begin{pmatrix}
  x_4\\
  x_5\\
  x_6
\end{pmatrix}
=0.
\]
Starting from $\overline{Y}_{2,3}$, we construct $Y$ as an
\emph{unprojection}, that is, the contraction of the Weil
divisor $\Pi \subset \overline{Y}$; see \cite{kinosaki},
\cite{MR1941620} for details.

\subsubsection{$h^0(\Delta, \l(2))=2$}
\label{sec:cb7-3}

Let $N \to \Delta$ be as in \ref{sec:case-3}. In this case, we can
immediately write down equations for a conic bundle with this
ramification data. Indeed, $A$ is a
$3\times 3$ symmetric matrix of homogeneous forms of degrees
  \[
  \begin{pmatrix}
    3&3&2\\
     &3&2\\
     & &1
  \end{pmatrix}.
  \]
We can take
  \[
  X= \{{}^t \mathbf{x} A \mathbf{x}=0 \} \in |2M+L| \subset \F(1,1,0).
  \]

We now exhibit an alternative model of $X$ as a $dP_3$ fibration.
We begin with a preliminary discussion of the geometry of the scroll 
$\F=\F(1,1,0)$, which we think of as a \emph{geometric} quotient
  \[
  \F= \a  /\!\!/ G
  \]
of $\a=\c^6$, with coordinates $u_0,u_1,u_2, x_0,x_1,x_2$, by the
group 
$G=\c^\times \times \c^\times$, with coordinates $\lambda,
\mu$, acting in the usual way (see the Appendix). Indeed, 
a $G$-linearisation of the trivial line bundle on $\a$ is the same
as a character $\chi\colon G \to \c^\times$. The group $G$ acts on
global sections $f\colon \a \to \c$ by $gf(x)=\chi(g)f(g^{-1}x)$.  
We say that a linearisation is \emph{useful} if the (open) subset 
 \[
\a^{\textit{ss}}_\chi =
\{x \in \a \mid \exists \; f \in \o^\chi_\a,\;f(x) \not = 0\}
 \]
of semistable points is nonempty, where 
\[
\o_\a^\chi = \{f \colon \a \to
\c \mid f(gx)=\chi(g) f(x)\}
\]
is the set of $G$-invariant sections. 
We have that
  \[
  \a^{\textit{ss}}_\chi=(\c^3 \setminus \{0\}) \times (\c^3 \setminus \{0\})
  \quad \Longleftrightarrow \quad \chi \in \r_+[L]+\r_+[M].
  \]
There are other useful linearisations. In fact, the cone of useful
linearisations is the cone $\r_+[L]+\r_+[M-L]$. This cone is
naturally partitioned in two chambers, and 
  \[
  \a^{\textit{ss}}_\chi=(\c^4 \setminus \{0\}) \times (\c^2 \setminus \{0\})
  \quad \Longleftrightarrow \quad \chi \in \r_+[M]+\r_+[M-L].
  \]
gives the second chamber.
When choosing a linearisation in the second chamber, the geometric
quotient is a different scroll over $\p^1$:
  \[
\a /\!\!/ G=\F^\prime = \F(0,1,1,1).
  \]
Crossing the wall separating the two chambers corresponds to a
birational map $\F \dasharrow \F^\prime$. The equation of
$X$, viewed as the equation of a hypersurface $Y \subset
\F^\prime$, is a section of the line bundle $\o_{\F^\prime} (3M -L)$.
The most convincing way to see this is to make the substitutions
  \begin{align*}
    u_0, u_1, u_2, x_0 & \mapsto y_1,y_2,y_3,y_0 \cr
    x_1,x_2 & \mapsto t_0,t_1
  \end{align*}
and to think of $t_0,t_1,y_0,y_1,y_2,y_3$ as natural coordinates on
the scroll $\F^\prime$; the action of $\lambda, \mu$ on these
coordinates is
\begin{align*}
 \lambda \colon (t_0,t_1,y_0,y_1,y_2,y_3) \mapsto & (\lambda^{-1}
 t_0,\lambda^{-1} t_1, y_0,\lambda y_1,\lambda y_2, \lambda y_3)\cr
 \mu \colon (t_0,t_1,y_0,\dots,y_3) \mapsto & (\mu t_0,\mu
 t_1,\mu y_0,y_1,y_2,y_3).
\end{align*}
To recover the standard presentation of the scroll $\F^\prime =
\F(0,1,1,1)$, we
change coordinates in $G$: the 1-parameter subgroups 
$\lambda^\prime \to (\lambda^{\prime -1},1) $ and 
$\mu^\prime \to (\mu^\prime,\mu^\prime)$ act
as 
\begin{align*}
 \lambda^\prime \colon (t_0,t_1, y_0,y_1,y_2,y_3) \mapsto & (\lambda^\prime
 t_0,\lambda^\prime t_1, y_0,\lambda^{\prime -1} y_1,\lambda^{\prime
 -1} y_2, \lambda^{\prime -1} y_3)\cr
 \mu^\prime \colon (t_0,t_1,y_0,y_1,y_2,y_3) \mapsto & (t_0, t_1,
 \mu^\prime y_0,\mu^\prime y_1,\mu^\prime y_2, \mu^\prime y_3).
\end{align*}
These calculations show that, with the stated substitutions, 
$H^0\bigl(\F^\prime ,3M-L\bigr)$ is canonically identified with
\[
H^0\bigl(\F, 3M-(-L+M)\bigr) = H^0(\F, 2M+L\bigr).\]
It is easy to see that the map $X \dasharrow Y$ is a flop which, by
what we just said,  is a Sarkisov link (of type IV according to 
\cite{MR96c:14013}) from the conic bundle
$X/\p^2$ to a $dP_3$ fibration $Y/\p^1$. 

\subsubsection{$h^0(\Delta, \l(2))=3$}
\label{sec:cb7-4}

Let $N \to \Delta$ be as in \ref{sec:case-4}. In this case, $A$ is a
$3\times 3$ symmetric matrix of homogeneous forms of degrees
  \[
  \begin{pmatrix}
    5&3&3\\
     &1&1\\
     & &1
  \end{pmatrix}
  \]
and we can take
  \[
  X= \{{}^t \mathbf{x} A \mathbf{x}=0 \}\in |2M+L| \subset \F(2,0,0).
  \]
The divisor $E=\{x_2=0\} \cap X$ is the exceptional divisor of 
a (2,1)-contraction $g\colon X \to Y$ with $K=0$, that is, $E$
contracts to a curve
of strictly canonical singularities on $Y$. This is a bad link and no
alternative model as a Mfs is produced.  We suspect that $X$ is
birationally rigid.  Observe that this follows from Conjecture~\ref{con:cb}.

\subparagraph{}

Table~\ref{tab:1} summarises the examples discussed so far. 
\begin{table}[ht]
\begin{center}
$
\begin{array}{|c|c|c|c|}
\hline
h^0\bigl(\Delta,\l(2)\bigr) & \text{model over} \; \p^2 & \text{link} &
\text{other model} \\
\hline
\hline
0 &X_{2M+L,M+L,M+L}          & \text{flop and} & Y_{2,2,2} \subset \p^6\\
  & \subset \p^2 \times \p^4 & \text{(2,1)-contraction}  & \\
\hline
1 &X_{2M+L,M+L}        &\text{flop and} & Y_{2,3,3,3,3}\subset \p(1^6,2) \\ 
  & \subset \F(1, 0^3) &\text{(2,0)-contraction} & \\
\hline
2 &X_{2M+L} \subset \F(1,1,0)& \text{flop}
& Y_{3M-L} \subset \F(1^3,0)\\
  & & & \\
\hline
3 &X_{2M+L} \subset \F(2,0,0)& \text{bad $K$ trivial}
&\text{$X$ rigid ?} \\
 & &\text{(2,1)-contraction} & \\
\hline
\end{array}
$
\end{center}
\caption{Conic bundles over $\p^2$ with discriminant $\Delta$ of
  degree 7\label{tab:1}}
\end{table}

\subsection{Conjectures on rigid conic bundles}
\label{sec:rigid-conic-bundles}

We briefly sketch a few more examples of conic bundles, mainly
nonstandard or over surfaces other than $\p^2$, then state some conjectures.

\begin{exa}
  We show that a general codimension 3 Fano \mbox{3-fold} 
  $Y_{3,3,4,4,4}\subset \p(1^5,2,3)$ is
  linked to a \emph{nonstandard} conic bundle $X$. 
  The equations of $Y$ are the Pfaffians of a $5\times 5$
  antisymmetric matrix of homogeneous forms which, in suitable
  coordinates, can be written as follows:
\[
\begin{pmatrix}
  z &a_1 &a_2 &a_3 \cr
    &b_1 &b_2 &b_3 \cr 
    &    &x_5 &x_4 \cr
    &    &    &x_3 
\end{pmatrix}.
\]
 (Here as usual $z$ is a coordinate of weight $3$, $y$ is a coordinate
 of degree $2$, and $x_1,\dots,x_5$ are coordinates of weight 1. The $a$s and
 the $b$s are homogeneous forms of degree 2 in the $x$s and $y$.)
 Projecting $Y$ to $\p(1^5, 2)$ we obtain a complete intersection
 $\overline{Y}_{3, 3}$ given by equations:
\[
\begin{pmatrix}
  a_1 &a_2 &a_3\cr
  b_1 &b_2 &b_3
\end{pmatrix}
\begin{pmatrix}
  x_3 \\
  x_4 \\
  x_5
\end{pmatrix}
=0,
\]
containing the weighted plane $\Pi =\{ x_3=x_4=x_5=0\}$. We now
construct a model as a conic bundle. First we construct a suitable
ambient space. Consider the quotient 
\[
\F = \bigl((\c^3 \setminus \{0\}) \times (\c^4 \setminus \{0\})\bigr)/
(\c^\times \times \c^\times)
\]
  by the action
\begin{align*}
  \lambda \colon (u_0,u_1,u_2,x_0,x_1,x_2,y) & \mapsto (\lambda u_0,
  \lambda u_1, \lambda u_2, \lambda^{-1} x_0, x_1, x_2, y) \cr
  \mu \colon (u_0,u_1,u_2,x_0,x_1,x_2,y) & \mapsto (u_0,u_1,u_2,\mu x_0, \mu
  x_1, \mu x_2, \mu^2 y). 
\end{align*}
  We can think of $\F$ as a scroll over $\p^2$ with fibre the weighted
  projective space $\p(1^3, 2)$. We define a morphism $f\colon \F \to
  \p(1^5,2)$ by setting:
\[
(u_0,u_1,u_2,x_0,x_1,x_2,y) \mapsto (x_1,x_2,u_0x_0,u_1x_0,u_2x_0,y).
\] 
  The proper preimage $X =f^{-1} \overline{Y}$
  is a general complete intersection of the form 
  $X=X_{2M+L,2M+L}\subset \F$. The
  ramification data $N \to \Delta$ is a \mbox{2-to-1} covering of a
  special plane octic with a triple point $P\in \Delta$, so $X$
  is not a standard conic bundle. It is easy to see that $X$ has an
  index 2 orbifold point over $P$. Because of the singular point, the
  threshold invariant $\tau(X/\p^2)$ \cite[\S 4.1]{corti:00} 
  is $5/2$, and not the expected $8/2$.
  Blowing up the point $P\in \p^2$ expresses $X$ as a conic bundle
  over $\F^1$ with discriminant of relative degree $5$, and this in
  turn has a model as a $dP_4$ fibration over $\p^1$.   
\end{exa}

\begin{exa}
  A general codimension 3 Fano \mbox{3-fold} $Y_{4,4,4,4,4} \subset \p(1^3,
  2^4)$ is linked to a nonstandard conic bundle 
  $X/\p^2$. The discriminant $\Delta$ is a
  plane curve of degree $9$ with four ordinary triple points, so $X$ is
  not a standard conic bundle. It is easy to see that $X$ has
  index 2 orbifold singularities over the singular points of
  $\Delta$. Please do your own calculations.
\end{exa}

\begin{exa}
  Takagi and Reid \cite{MR1924722}, \cite{kinosaki} construct a
  codimension 4 Fano \mbox{3-fold} $Y \subset \p(1^6,2^2)$ with $h^0(-K)= g+2
  =6$ and $-K^3=7$ with 2 index 2 orbifold points. There are two
  deformation families with these invariants, corresponding to the
  ``Tom'' and ``Jerry'' formats of unprojection. If $Y$ is a general
  Fano in the Tom family, then $Y$ is linked to a general conic bundle
  $X\to \p^2$ with discriminant curve of degree 6. The cover $N \to
  \Delta$ corresponds to the partition $6=2+2+2$, but we did not
  yet carry out all the necessary calculations to understand the model
  $X$ explicitly.    
\end{exa}

\begin{exa}
  There is no reason why we should only work with conic bundles over
  $\p^2$, since Catanese's theory works essentially unchanged for curves on
  a weighted projective plane. For example, there are two
  Fano \mbox{3-folds} in codimension 4 and 5, $Y \subset \p(1^5,2^3)$ and
  $Y\subset \p(1^5,2^4)$, with $-K^3=11/2$ and $6$, see \cite{MR1924722}, 
  which can be linked to conic bundles over $\p(1,1,2)$ with discriminants in
  $\o(8)$ and $\o(6)$.
\end{exa}

We have written down explicit equations of several families of conic bundles
with discriminant of low degree. We believe that the \emph{explicit
geometry} of these varieties will play an increasingly more prominent
role in the study of their birational geometry. We state the following
rather optimistic conjecture.

\begin{dfn}
 We say that a standard conic bundle $\pi \colon X\to S$ satisfies 
condition~$(\ast)$ if $-K_X \not \in \Ins \NM^1 X$.   
\end{dfn}

\begin{con} \label{con:cb}
 A standard conic bundle over $\p^2$ is
 birationally rigid if it satisfies condition~$(\ast)$. 
\end{con}
We feel like making the following conjecture for which we have little
evidence. 
\begin{con}
 A standard conic bundle over $\p^2$ is
 birationally rigid if the discriminant has degree \mbox{$\ge 9$}.
\end{con}

\section{$dP_3$ fibrations }
\label{sec:dp_3-fibrations}

\begin{conv} \label{conv:dp3} 
In this section a \emph{$dP_3$ fibration} is a \mbox{3-fold} $X$
together with a 
morphism $f\colon X \to \p^1$ satisfying the following conditions:
\begin{enumerate}
\item The nonsingular fibres of $f$ are cubic surfaces.
\item $X$ has Gorenstein terminal singularities (these are precisely
  the isolated hypersurface singularities with a DuVal section). In particular
  $X\subset \F= \p (E)$ is naturally embedded in a rational scroll
  over $\p^1$ (a natural choice is $E=f_\ast \o (-K_X)$). 
\item $X$ has Picard rank $\rho = 2$, that is, the morphism $f\colon X
  \to \p^1$ is extremal.
\item  The local rings of $X$ are unique factorisation domains, 
  that is, Weil divisors on $X$ are Cartier.
\end{enumerate}
\end{conv}
\noindent In addition, we often assume that $X$ is nonsingular, or at
  least that $f\colon X \to \p^1$ is \emph{semistable} in the sense of
  \cite{MR98g:11076}. 

\begin{rem}
Corti and Koll\'{a}r \cite{MR98e:14037}
\cite{MR98g:11076} show that if $f\colon X\to \p^1$ is a Mori fibre space
$dP_3$ fibration, then $X/\p^1$ is square birational to a semistable 
$g\colon Y \to \p^1$; in particular $Y$ is Gorenstein and the
conditions above are satisfied. It is therefore not restrictive to
limit our attention to semistable fibrations.
\end{rem}

In this section, we aim to do two things. First, we want to determine
the geography of $dP_3$ fibrations, that is, determine all integers
$n,a,b,c$ such that a general member $X \in |3M+nL|$ on $\F(0,a,b,c)$
is a $dP_3$ fibration in our sense. 
Second, we want to state some
conjectures on the birational geometry of $dP_3$ fibrations.
M.~Grinenko has been studying del Pezzo
fibrations systematically in a series of recent papers
\cite{grinenko:01a, grinenko:01b, grinenko:01c}. Because he is
primarily concerned with fibres of degree 1 and 2, there is little
overlap between his work and what we do here. 

\subsection{The $K^2$~condition}
\label{sec:k2-cond}

\begin{dfn} \label{dfn:k2-cond}
  We say that $X$ satisfies the \emph{$K^2$~condition} if
\[ 
K_X^2 \not \in \Ins \NE X.
\] 
\end{dfn}

\begin{thm} \cite{MR98j:14014} \label{thm:dp3}
  Let $f\colon X\to \p^1$ be a $dP_3$ fibration
  satisfying the following technical conditions:
  \begin{itemize}
  \item the total space $X$ is nonsingular, and $f$ has Lefschetz
  singularities, that is, it has only ordinary critical points,
  with distinct critical values. 
  \item If $X_b$ is a singular fibre, then there are exactly six lines
  of $X_b$ passing through the unique singular point.   
  \end{itemize}
  If in addition $X$ satisfies the
  $K^2$~condition, then $X$ is birationally rigid.
\end{thm}

\begin{rem}
  It follows from the above theorem that, under its assumptions,
  Iskovskikh's conjecture \cite{MR2000i:14019} holds.
\end{rem}

If a $dP_3$ fibration $X\to S$, with $X$ nonsingular, 
is nonrigid, then $X\to S$ belongs to one of finitely
many algebraic families. (This can be proved using 
Proposition~\ref{pro:finite}.) On the other hand, as we see in
Section~\ref{sec:geog}, infinitely many families do not
satisfy the $K^2$~condition. This shows that Pukhlikov's theorem is
not optimal.

\subsection{Conjectures}
\label{sec:dp-con}

\begin{dfn}
  We say that $X$ satisfies the \emph{condition~$(\ast)$} if $-K \not \in \Ins
\NM^1 X$.
\end{dfn}

\begin{rem}
  In Section~\ref{sec:nonrigid} we study (among other things)
  condition~$(\ast)$ for the general members of families of
  $dP_3$. Though we don't prove it completely, we believe that
  condition~$(\ast)$ is satisfied by a general member of all but a
  handful of families of $dP_3$ fibrations listed in Table~\ref{tab:2}.
\end{rem}

Grinenko has recently made the following striking conjecture. 

\begin{con} \cite[Conjecture 1.5]{grin99} \cite[Conjecture 1.6]{grinenko:00} 
  \cite[Conjecture 2.5]{gri02} \label{con:dp-ast}
  A $dP_3$ fibration with nonsingular total space
  is birationally rigid if it satisfies condition~$(\ast)$.
\end{con}

In the remaining part of this subsection we make a few comments on the
meaning of the conjecture. 

\begin{con} \label{con:dp-ast-2}
  Same as \ref{con:dp-ast}, only assuming that $X/\p^1$ is
  semistable in the sense of \cite{MR98g:11076}. 
\end{con}

\begin{con} \label{con:dp-equiv}
  Let $X\to \p^1$ be a $dP_3$ fibration, with $X$ nonsingular. Let
  $X^\prime \to \p^1$ be a \mbox{3-fold} Mori fibre space, square birational
  to $X/\p^1$. If $X/\p^1$ satisfies condition~$(\ast)$, then so does 
  $X^\prime/\p^1$. 
\end{con}

\begin{con} \label{con:dp-equiv-2}
  Same as \ref{con:dp-equiv}, only assuming that $X/\p^1$ is
  semistable.
\end{con}
\noindent More experimentation is needed before we can have any confidence
  in these conjectures. Here we only briefly touch on these matters in
  the Example in Section~\ref{sec:unstable}. 

\begin{pro}
  Conjecture~\ref{con:dp-ast} follows from
  Conjecture~\ref{con:dp-equiv}, and Conjecture~\ref{con:dp-ast-2}
  from Conjecture~\ref{con:dp-equiv-2}. 
\end{pro}

\noindent This indicates that further progress is likely to
come from a systematic study of square birational maps involving a
semistable $dP_3$ fibration and a Mori fibre space. 

\begin{proof}
  We sketch the proof. If $X\to \p^1$ is not rigid, there is
  a Mfs $Y \to T$ and a nonsquare birational map $X \dasharrow
  Y$. Applying the Sarkisov program as in \cite{MR96c:14013} gives a
  Mfs $f^\prime \colon X^\prime \to \p^1$ square birational to 
  $X\to \p^1$, and a linear system 
\[
\mathcal{H}^\prime \subset |-nK_{X^\prime} + f^{\prime \ast} A|
\]
  where $A$ is a divisor on $\p^1$ of strictly negative degree. It follows that
  $X^\prime$ does not satisfy condition~$(\ast)$.
 \end{proof}

\subsection{Geography for $dP_3$ fibrations}
\label{sec:geog}

\subsubsection{Notation and basic numerology}
\label{sec:notat}

\paragraph{Notation}
A $dP_3$ fibration $X\to\p^1$, as defined formally in Convention~\ref{conv:dp3},
always admits model as a hypersurface in a 4-fold scroll $\F$
that is a $\p^3$ bundle over $\p^1$. We fix the notation in use
throughout this section.
\begin{enumerate}
\item $X \in |3M + nL|\subset \F = \F(0,a,b,c)$, where $a,b,c,n$ are
 integers with  $0\le a \le b \le c$. We write $d=a+b+c$.  
\item We write $u,v,x,y,z,t$ for the homogeneous coordinates of $\F$,
  where $u,v$ are the homogeneous coordinates on the base $\p^1$ and
  $x,y,z,t$ are the fibre coordinates.
\item We denote $L$ and $M$ the natural basis of $\Pic (\F)$. These
  line bundles have sections $u\in H^0(\F, L)$ and $x\in H^0(\F, M)$,
  and we sometimes identify $L$, $M$ with the actual divisors $u=0$,
  $x=0$.
\item We write $\Gamma = \{y=z=t=0\}$ and $B=\{z=t=0\}$. Note that
  $\Gamma$ generates an extremal ray of $\NE (\F)$. 
\item We denote by $F\in \o_\F (3M+nL)$ the equation of $X$. The
  polynomial $F$ is a sum $\sum\alpha m$ of \emph{terms} $\alpha m$
  where the sum ranges over the \emph{fibre monomials} $m$ that are
  cubics in $x,y,z,t$ with coefficients $\alpha=\alpha(u,v)$ to
  fix up homogeneity.  We write both $\alpha m\in F$ and $m\in F$ to
  mean that the term (with implicit coefficient $\alpha$ if not
  mentioned) appears in $F$ with nonzero coefficient.
\end{enumerate}

\paragraph{Basic numerology}
\begin{enumerate}
\item $M^3L = 1$, $M^4 = d$.
\item $-K_X = -K_{|X}$, where $-K = M + (2 - d - n)L$ is a divisor on $\F$.
\item $X \cdot \Gamma = n$. Apart from the trivial case $X\in |3M|$ on
  $\p^1 \times \p^3$ (a constant family of cubic surfaces) there are
  two cases:
\[
  \begin{cases}
    n\geq 0 &\quad \text{and $3M+nL$ is nef and big},\\
    n <   0. &
  \end{cases}
\]
 When $3M+nL$ is nef and big, it is base point free and a general $X
 \in |3M+nL|$ is nonsingular. It then follows from the Lefschetz
 hyperplane theorem that $\rho (X)=2$ and $X$ is a $dP_3$ fibration in
 our sense. Almost everything we say below refers to the much more
 interesting case when $n<0$.
\end{enumerate}

\subsubsection{What the picture says} \label{sec:geog_pic}

\begin{figure}
\setlength{\unitlength}{1.2mm}
\begin{center}
\begin{picture}(100,155)
\put(0,5){\vector(1,0){95}}
\put(65,0){\vector(0,1){150}}
\put(95,0){\makebox(5,10){$n$}}
\put(60,150){\makebox(10,5){$d$}}
\put(82,0){\makebox(5,5){$2$}}
\put(72,0){\makebox(5,5){$1$}}
\put(64,0){\makebox(5,5){$0$}}
\put(52,0){\makebox(5,5){$-1$}}
\put(42,0){\makebox(5,5){$-2$}}
\put(32,0){\makebox(5,5){$-3$}}
\put(22,0){\makebox(5,5){$-4$}}
\put(12,0){\makebox(5,5){$-5$}}
\put(2,0){\makebox(5,5){$-6$}}
\put(65,40){\makebox(5,10){$4$}}
\put(65,60){\makebox(5,10){$6$}}
\put(65,80){\makebox(5,10){$8$}}
\put(65,100){\makebox(5,10){$10$}}
\put(65,120){\makebox(5,10){$12$}}
\put(65,140){\makebox(5,10){$14$}}
\put(65,45){\line(3,-5){25}}
\put(65,45){\line(-3,5){65}}
\put(65,25){\line(1,-1){25}}
\put(65,25){\line(-1,1){20}}
\put(45,45){\line(-1,2){10}}
\put(35,65){\line(0,1){10}}
\put(35,75){\line(-1,1){20}}
\put(15,95){\line(-1,2){10}}
\put(5,115){\line(0,1){10}}
\put(5,125){\line(-1,1){5}}
\multiput(84.5,4.2)(0,10){15}{$\cdot$}
\multiput(74.5,34.2)(0,10){12}{$\cdot$}
\multiput(64.5,54.2)(0,10){10}{$\cdot$}
\multiput(54.5,64.2)(0,10){9}{$\cdot$}
\multiput(44.5,84.2)(0,10){7}{$\cdot$}
\multiput(34.5,94.2)(0,10){6}{$\cdot$}
\multiput(24.5,114.2)(0,10){4}{$\cdot$}
\multiput(14.5,134.2)(0,10){2}{$\cdot$}
\put(74.3,4){$\bullet$}
\put(74.3,14){$\bullet$}
\put(64.3,14){$\bullet$}
\put(64.3,24){$\bullet$}
\put(54.3,34){$\bullet$}
\put(54.3,44){$\bullet$}
\put(44.3,44){$\bullet$}
\put(44.3,54){$\bullet$}
\put(44.3,64){$\bullet$}
\put(34.3,64){$\bullet$}
\put(34.3,74){$\bullet$}
\put(74.3,24){$\circ$}
\put(64.3,34){$\circ$}
\put(64.3,44){$\circ$}
\put(54.3,54){$\circ$}
\put(44.3,74){$\circ$}
\put(34.3,84){$\circ$}
\put(34.3,94){$\circ$}
\put(24.3,84){$\circ$}
\put(24.3,94){$\circ$}
\put(24.3,104){$\circ$}
\put(14.3,94){$\circ$}
\put(14.3,104){$\circ$}
\put(14.3,114){$\circ$}
\put(14.3,124){$\circ$}
\put(4.3,114){$\circ$}
\put(4.3,124){$\circ$}
\put(4.3,134){$\circ$}
\put(4.3,144){$\circ$}
\put(69,10){\makebox(10,5){{\footnotesize $001$}}}
\put(57,21){\makebox(10,5){{\footnotesize $011$}}}
\put(49,30){\makebox(10,5){{\footnotesize $111$}}}
\put(50,40){\makebox(10,5){{\footnotesize $[112]$}}}
\put(42,47){\makebox(10,5){{\footnotesize $[113]$}}}
\put(39,40){\makebox(10,5){{\footnotesize $112$}}}
\put(40,50){\makebox(10,5){{\footnotesize $122$}}}
\put(28,60){\makebox(10,5){{\footnotesize $123$}}}
\put(40,60){\makebox(10,5){{\footnotesize $[123]$}}}
\put(27,71){\makebox(10,5){{\footnotesize $133$}}}
\put(49,54){\makebox(10,5){{\footnotesize $(122)$}}}
\put(41,64){\makebox(10,5){{\footnotesize $(222)$}}}
\put(44,75){\makebox(10,5){{\footnotesize $(133)$}}}
\put(32,75){\makebox(10,5){{\footnotesize $(223)$}}}
\put(31,84){\makebox(10,5){{\footnotesize $(233)$}}}
\put(33,94){\makebox(10,5){{\footnotesize $(144)$}}}
\put(23,84){\makebox(10,5){{\footnotesize $(224)$}}}
\put(21,94){\makebox(10,5){{\footnotesize $(234)$}}}
\put(19,104){\makebox(10,5){{\footnotesize $(244)$}}}
\put(13,94){\makebox(10,5){{\footnotesize $(225)$}}}
\put(11,104){\makebox(10,5){{\footnotesize $(235)$}}}
\put(11,114){\makebox(10,5){{\footnotesize $(245)$}}}
\put(15,124){\makebox(10,5){{\footnotesize $(255)$}}}
\put(4,114){\makebox(10,5){{\footnotesize $(236)$}}}
\put(3,124){\makebox(10,5){{\footnotesize $(246)$}}}
\put(1,134){\makebox(10,5){{\footnotesize $(256)$}}}
\put(4,143){\makebox(10,5){{\footnotesize $(266)$}}}
\put(3,150){\parbox[b]{70mm}{The Pukhlikov line $3d + 5n = 12$}}
\put(0,55){\parbox[b]{30mm}{Symbols:}}
\put(0,50){\parbox[b]{50mm}{$\cdot$ possible $dP_3$ coordinate}}
\put(0,45){\parbox[b]{50mm}{$\circ$ for $X$ with a $K=0$ bad link}}
\put(0,40){\parbox[b]{50mm}{$\bullet$ for known nonrigid $X$}}
\put(0,30){\parbox[b]{30mm}{Labels:}}
\put(0,25){\parbox[b]{70mm}{$abc$ when every $X\subset\F(a,b,c)$ admits the link}}
\put(0,20){\parbox[b]{70mm}{$[abc]$ when only special $X$ admit the link}}
\put(0,15){\parbox[b]{70mm}{$(abc)$ when every $X$ has $K=0$ bad link}}
\end{picture}
\end{center}
\caption{Geography of $|3M+nL|\subset\F(a,b,c)$ with $d=a+b+c$ 
\label{fig:geog}}
\end{figure}
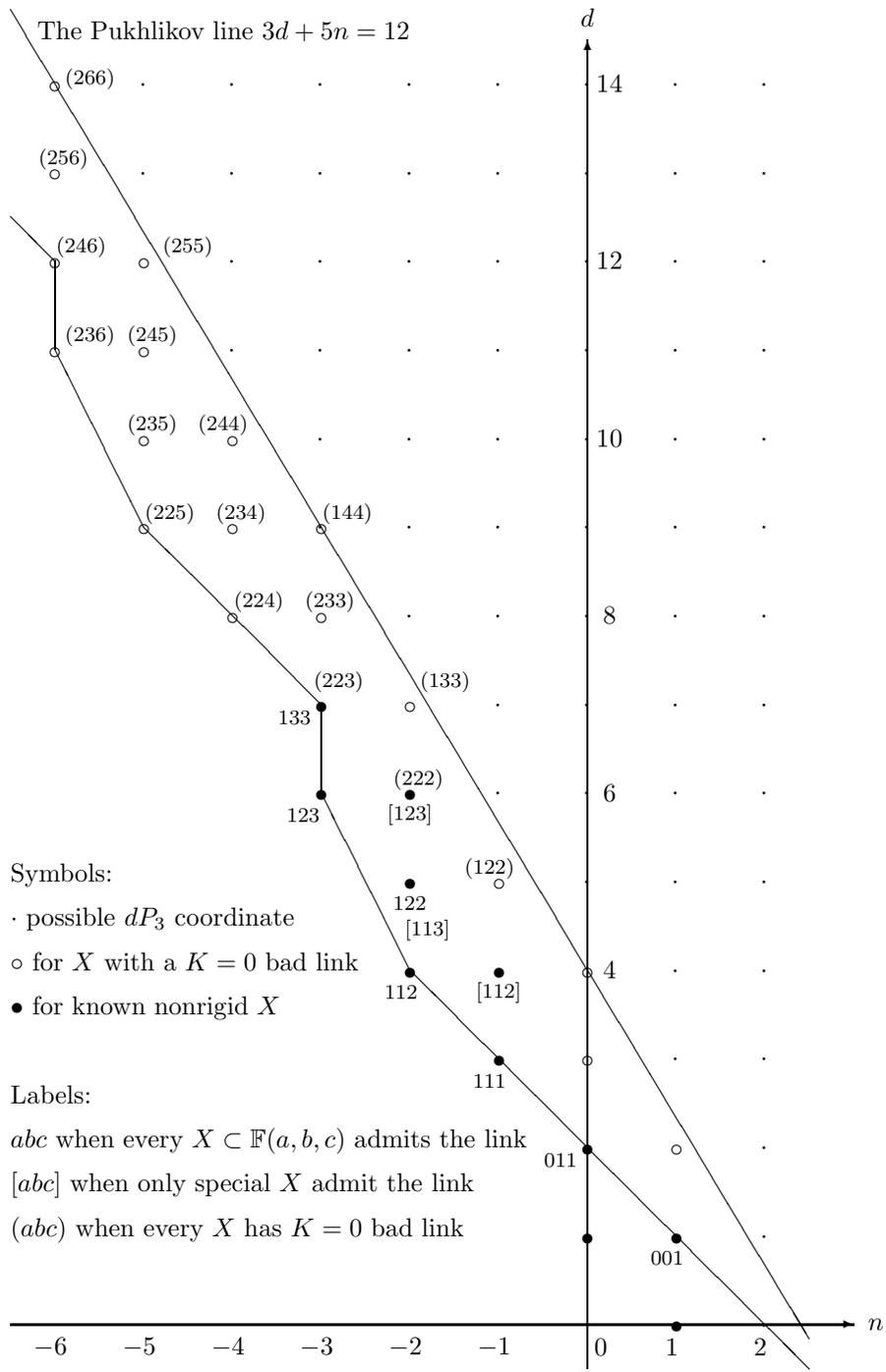

We plot the families of $dP_3$ fibrations on 
a graph of $n$ against $d$ as in Figure~\ref{fig:geog} and we refer to this
picture as \emph{geography}. We regard the figure as a pictorial
statement of a theorem that we now spell out.

The figure summarises information
obtained in calculations carried out in the remaining part of the
paper; in particular the figure displays the following.
\begin{enumerate}
\item Pairs $(n,d)$ for which there is a family of $dP_3$ fibrations $X\in
  |3M+nL|$ in $\F(0,a,b,c)$ with $d=a+b+c$.
\item Triples $(a,b,c)$ such that a general or a special member of the
  corresponding family is known to be nonrigid.
\item Triples $(a,b,c)$ such that a general member of the
  corresponding family does not satisfy condition~$(\ast)$. We
  do not prove it completely, but we believe that the figure displays
  all such triples. 
\end{enumerate}

We plan to use the picture as a primary testing ground for 
Grinenko's Conjecture and as a starting point possibly to prove
it. Here we explain what the picture and its
various elements mean. Precise details are worked out in the following
subsections.

\subparagraph{Families of $dP_3$ fibrations}

We mark points $(n,d)$ of Figure~\ref{fig:geog} by a dot $\cdot$
if and only if there are values $0\leq a \leq b \leq c $ with
$d=a+b+c$, such that a general $X \in |3M+dL|$ in $\F(0,a,b,c)$ is a 
$dP_3$ fibration in our sense. A dot can be a bullet~$\bullet$ or a
circle~$\circ$; we explain what these mean below. The geography
consists of the first quadrant minus the origin, and the region to the
right of a curve with a periodic behaviour which we define by the
picture itself. A point of the picture may house several deformation families: 
for example, $(n,d)$ is $(-2,5)$ for a family of $dP_3$ fibrations in
$\F(1,1,3)$ and also for a different family that lies in $\F(1,2,2)$.
We explain how this region is drawn in
Section~\ref{sec:how-geogr-obta} below.

\subparagraph{Nonrigid families}

A bullet~$\bullet$ marks a point $(n,d)$ corresponding to
some family of $dP_3$ fibrations for which we know that a member $X$ is 
nonrigid. When $X$ is general in its family, we specify the family 
by writing $abc$ under the bullet, indicating the
invariants of $\F(0,a,b,c)$. In some cases, we only know that a 
\emph{special} $X$ in the family is nonrigid. We indicate this by
writing $[abc]$ instead. In all cases when we know that $X$ is
nonrigid, the 2-ray game from $X/\p^1$ is a nonsquare link of the
Sarkisov program: we refer to this as \emph{the} link
and we say that \emph{$X$ admits the link}. We describe the link
explicitly in Section~\ref{sec:nonrigid} and summarise our calculations in 
Table~\ref{tab:2}.

\subparagraph{Condition~$(\ast)$}

We do not prove it completely, but we believe that
  condition~$(\ast)$ does not hold for a general member 
  $X\in |3M+nL|$ in $\F(a,b,c)$
  if and only if $(n,d)$ has $\bullet$ with $abc$. 
  After the argument of Section~\ref{sec:mfs_con}, we only study points
  to the left of the line $3d + 5n - 12 = 0$; the arguments of
  Section~\ref{sec:family5} show immediately that about two thirds
  of those with $n<3$ satisfy condition~$(\ast)$, and we believe
  the remaining cases to be only slightly more difficult.
  Grinenko's
  conjecture is then equivalent to the statement that Table~\ref{tab:2}
  is a complete list of families of $dP_3$
  fibrations with nonrigid general element.  
  
\subparagraph{$K$ trivial contractions}

  Unless there is already a bullet at this point, we mark $(n,d)$
  with a circle $\circ$ when we know that 
  there is a corresponding family of $dP_3$ 
  fibrations such that the 2-ray game starting with a general member 
  $X$ terminates with a $K$-trivial contraction and a variety with
  strictly canonical singularities. We specify the family by writing 
  $(abc)$ above the circle.  Note that, in each example
  of $[abc]$, a general $X$ in the family has a $K$-trivial bad link
  and we do not repeat $(abc)$ there. These examples are not our main
  interest here, but they are cases where the question of
  rigidity is particularly intriguing. 
 
\subparagraph{The Pukhlikov line}

  The figure also shows the line $3d+5n=12$. As we show in
  Lemma~\ref{lem:explicitk2} below, the
  general member of a family satisfies the $K^2$~condition 
  if the point $(n,d)$ lies to the right of this line. (For $n$
  negative, this is  if and only if, which leaves out a couple of
  cases with $n>0$ where we are not completely sure.) The graph shows
  that the line leaves out a thin strip of the geography, containing
  infinitely many families. 

\subsubsection{How the geography is obtained}
\label{sec:how-geogr-obta}

In this subsection we make the statements from
which the geography of $dP_3$ fibrations is derived.
\begin{pro} \label{pro:geog}
  Figure~\ref{fig:geog} marks all pairs $(n,d)$ for which there are integers
  $0\leq a \leq b \leq c$ such that: 
  \begin{enumerate}
  \item The relative surface $B=\{z=t=0\}$ is not contained in the
  base locus of the linear system $|3M+nL|$. (When $a=b$, we require
  that a general member of the linear system does not contain a
  surface of the form $\{s=t=0\}$ where $s$ is some section of $\o_\F (M-bL)$.)
  \item The curve $\Gamma = \{y=z=t=0\}$ is contained in the base
  locus of $|3M+nL|$ with multiplicity at most one, that is, when
  $n<0$, a general member $X$ in the linear system is nonsingular
  generically along $\Gamma$. 
 \end{enumerate}
\end{pro}
\noindent The conditions imply that a general member $X\in |3M+nL|$ is
nonsingular outside $\Gamma$ and possibly has isolated singularities 
along $\Gamma$. In fact, $X$ is a $dP_3$ fibration.
\begin{pro} \label{pro:dp3}
  Let $(n,d)$ and $X \in |3M+nL|$ in $\F(0,a,b,c)$ be as in 
  Proposition~\ref{pro:geog}
  above. Then:
  \begin{enumerate}
  \item For any $X$, $\Pic (X)= H^2(X) = \z^2$, that is, $X$ has
  Picard rank 2. 
  \item If all fibres of $f\colon X \to \p^1$ are reduced and
  irreducible, we have an exact sequence
\[
0\to \z [L] \to A_2 (X) \to \Pic (X_\eta) \to 0 
\]
  where $A_2(X)$ is the group of 2-dimensional cycles, that is, Weil
  divisors on $X$ modulo rational equivalence, and $X_\eta$ is the
  generic fibre. 
  \item If $X$ is general, then $\Pic (X_\eta) = \z[-K]$.
  \item If $X$ is general, then $X$ has isolated singularities of $cA$
  type.  
  \end{enumerate}
\end{pro}
\begin{cor}
  With the same assumptions, if $X$ is general, then $X/\p^1$ is a
  $dP_3$ fibration. 
\end{cor}

\paragraph{Proof of Proposition~\ref{pro:geog}}

The proposition is a direct consequence of the following more precise lemma.
Indeed, it is elementary to check that the inequalities in the lemma specify
the region marked in Figure~\ref{fig:geog}.

\begin{lem} \label{lem:ineqs}
  Consider the linear system $|3M+nL|$ on $\F(0,a,b,c)$ as above. 
  \begin{enumerate}
  \item $B \not \subset \Bs |3M+nL|$ if and only if $n\geq -3a$. 
  \item If $n<0$ and $a=b$, then all $X$ contain a surface
    $t=\ell(z,y)=0$, where $\ell$ is a linear form if and only if $n=-3a$.
  \item If $n<0$, so that $\Gamma \subset \Bs|3M+nL|$, then a
  general member $X\in |3M+nL|$ is nonsingular generically along
  $\Gamma$ if and only if $n\geq c$. In this case also $n+d \geq 2a$.
 \end{enumerate}
\end{lem}

\begin{proof}
If the first inequality fails, then every monomial in $F$ is divisible
by $z$ or $t$ so $X$ contains $B$.  

If $-n=3a>0$ and $a=b$, we can refine this analysis slightly. 
We can write $F=f_3(y,z)+tf_2(x,y,z,t)$, where $f_2$, $f_3$ are
homogeneous forms of the indicated degrees. If $\ell(y,z)$ is a linear
factor of $f_3$, then $X$ necessarily contains the surface $t=\ell(z,y)=0$.

Suppose that $n<0$, so that $a>-(1/3)n>0$ by the first inequality.
The fibre monomial $x^3$ cannot now occur in any
term of $F$, so the inclusion $\Gamma\subset X$ is clear.

The total space $X$ is nonsingular generically along $\Gamma$ if and
only if, for general values of $u,v$, the polynomial $F$ contains
at least one of the terms $x^2y$, $x^2z$, $x^2t$.  
If $x^2t\in F$, then $n \ge -c$,
and that is the weakest inequality of the three.  Finally
\[
d + n = (a+b) + (c+n) \ge 2a
\]
with equality only if $n=-c$ and $a=b$.
\end{proof}

\begin{rem}
Some of our statements in Section~\ref{sec:geog_pic} follow
from Lemma~\ref{lem:ineqs}.
Let $\sigma$ be the closed subcone of $\NM^1(X)$ generated
by the mobile rays $L$, $D_z=(z=0)\cap X$.  Then writing
$-K_X = (2-a-c-n)L + D_z$ shows that
\begin{itemize}
\item[(a)]
$2-a-c-n > 0
\quad
\text{if and only if}
\quad
-K_X \in\Ins\sigma
$,
\item[(b)]
$2-a-c-n \ge 0
\quad
\text{if and only if}
\quad
-K_X \in\sigma
$.
\item[(c)]
If condition~$(\ast)$ holds for $(n;a,b,c)$, then $2-a < c+n$.
\end{itemize}
It is easy to solve these inequalities for $a,c,n$ with the results
of Lemma~\ref{lem:ineqs}.  In case (a), the solutions are exactly those
4-tuples $(n;a,b,c)$ for which $(n,d)$ is marked by $\bullet_{abc}$.
In case (b), the solutions are exactly those $(n;a,b,c)$ for which
$(n,d)$ is marked by $\circ^{(abc)}$.  Condition~$(\ast)$ holds for
each of the solutions that have $a=2$; it fails only at the special
member of some solutions when $a=1$, which are the [abc] cases.
Part (c), which is virtually no condition at all when $n\ge 0$.
\end{rem}

\paragraph{Proof of Proposition~\ref{pro:dp3}}

We briefly sketch the proof, which is standard.
We take on the first assertion first. If $n\geq 0$, then $X$ is nef
and big. Denote $U=\F \setminus X$ the open complement. 
The Lefschetz hyperplane theorem implies that
$H^{8-i}_c(U)=H_i(U)=(0)$ for $i>4$ hence $H^2_c(U) = H^3_c (U)=(0)$
and from the standard exact sequence
\[
H^2_c(U) \to H^2 (\F) \to H^2 (X) \to H^3_c(U)
\]
we deduce that $H^2(X)=H^2(\F)=\z^2$. Since $H^1(X, \o_X)=(0)$, it
follows that $\Pic(X)=\z^2$, based by $L$ and $M$. When $n<0$
$X$ is not nef and we can proceed in various ways. For example, we can
observe that the proper transform $X^\prime \subset \F^\prime$ is nef
and big on the variety $\F^\prime$ obtained from $\F$ by flipping the
curve $\Gamma$ (see the Appendix for information on $\F^\prime$). By
Lefschetz then $H^2(X^\prime)= \z^2$. The birational map 
$X^\prime \dasharrow X$ is either an isomorphism in codimension 1, or
else $X^\prime$ contains the whole of the flipped set $\p(a,b,c)$. In
either case, we conclude that $H^2 (X)=\z^2$ and the statement follows.

The exact sequence in 2.\ is standard.

Consider now the generic fibre $S$ of $X\to\p^1$ which is a nonsingular
cubic surface embedded in $\p^3$ over the function field $\c(s)$
of $\p^1$ with coordinates $x,y,z,t$.  See \cite{KSC}
for a tutorial in the elementary methods of such cubics.  We claim
that $S$ has Picard rank~1.  An example of such $S$ is the surface
\[
S_0 = py^3 + qz^3 + rt^3 + x^2t
\]
where $p,q,r\in \c(s)$ are chosen generally.
It is easy to calculate the 27 lines on such a cubic surface,
and then to figure out the action of the Galois group $\Gal(K/\c(s))$
for any extension $K/\c(s)$ in which the lines are geometric:
for general $p,q,r$, the Galois group is $\z/3\times \z/3 \times \z/2$ and
the 27 lines split into orbits of size 3,6,9,9.  One checks
by explicit computation that the smaller orbits are made of
unions of entire hyperplane sections over $K$ so there are no
Galois orbits of disjoint lines.  Segre's theorem, \cite{KSC}
Theorem~2.16, then implies that the Picard rank of $S_0$ is 1.
Since a special surface $S_0$ in the family has rank~1, so
does the general surface $S$. This shows statement 3.

The final claim follows from an elementary monomial argument. Indeed $x^2t\in F$
  with coefficient a homogeneous polynomial $\varphi (u, v)$. The
  singularities of $X$ along $\Gamma$ are at the zeros of $\varphi$
  and if $X$ is general we may assume that these all have
  multiplicity 1. If $u|\varphi$, for example, the singularity at
  $u=0$ is of the form $ut+\cdots =0$. This implies that the general
  surface section has a singularity of type $A_k$ there.  \qed

\subsubsection{The $K_X^2$ condition and Pukhlikov's Theorem}

\begin{lem}
\label{lem:explicitk2}
If $X\to\p^1$ is a $dP_3$ fibration, then
\[
K_X^2\in\Ins\NE(X)
\quad
\text{implies that}
\quad
3d + 5n < 12.
\]
If in addition $n<0$, then
\begin{itemize}
\item[(a)]
the 1-cycles $M^2L$ and $\Gamma\subset X$ bound
the Mori cone $\NE(X)$
\item[(b)]
$K_X^2\in\Ins\NE(X)$ if and only if $3d + 5n < 12$.
\end{itemize}
\end{lem}
\begin{proof}
Suppose first that $n<0$.
The general fibre of $X\to\p^1$ is a cubic surface, so $M^2L$,
a line in the fibre, lies on one extreme ray of $\NE(X)$.
Since $a>0$, there is a projective morphism from the scroll $\F$
that contracts $\Gamma$, so the curve $\Gamma\subset X$ generates
the other ray of $\NE(X)$, proving (a).  

Part (b) follows from (a) by calculating the 1-cycle $K_X^2$ in the
basis $\Gamma, M^2L$.  Denoting $i\colon X \hookrightarrow\F$ for the
inclusion and noting that
\[
\Gamma=(M-aL)(M-bL)(M-cL)=M^3 - dM^2L,
\]
we calculate $i_\ast K_X^2=(-K)^2X$ on $\F$ as
\[
\begin{array}{rcl}
\bigl(M + (2 - d - n)L\bigr)^2(3M + nL) & = & 3M^3 + \bigl(6(2-d-n) + 
n\bigr)M^2L   \\
 & = & 3\Gamma + (12 - 3d - 5n)M^2L.
\end{array}
\]

Now let $n$ be any integer.  Suppose that $K_X^2\in\Ins\NE(X)$.  We
calculate in ${\mathrm N}_1(\F)$ omitting $i_\ast$.  Let $\sigma$ be the
closed subcone of effective divisors spanned by $K_X^2$ and $M^2L$.
If $\Gamma\notin\sigma$, then $K_X^2$ is a strictly convex combination
of $\Gamma$ and $M^2L$ so the inequality holds.  If $\Gamma\in\sigma$,
then some multiple $\Sigma=k\Gamma$ of the class of $\Gamma$ (on $\F$)
contains an effective curve that lies on $X$.  But $\Gamma\subset\F$
can be contracted, so $\Sigma$ must be supported on $\Gamma$.
Therefore $\Gamma\subset X$ and $\Gamma$ generates an extremal ray
of $\NE(X)$. The inequality follows as before.
\end{proof}

\subsection{Nonrigid $dP_3$ fibrations}
\label{sec:nonrigid}

We summarise all the examples of nonrigid $dP_3$ fibrations that
we know in Table~\ref{tab:2}.  We do not dwell on the well-known cases:
\begin{itemize}
\item
$(n,d)=(0,1)$ links to the cubic \mbox{3-fold} $Y=Y_3\subset\p^4$,
the $dP_3$ fibration being the pencil of $\p^3$s through a plane
intersecting $Y$ in a cubic curve
\item
$(n,d)=(1,0)$ links to $Y=\p^3$, the $dP_3$ fibration being the
pencil of any pair of transverse cubic surfaces in $Y$
\item
$(n,d) = (1,1)$ links to $Y_{3,3}\subset\p(1^5,2)$, the $dP_3$
fibration being a pencil of divisors having maximal vanishing
at the singular point; see \cite{cm02}, \cite{BZ03}.
\end{itemize}
We discuss the table below and then make three detailed studies
of examples.

\begin{table}[ht]
\begin{center}
\renewcommand{\arraycolsep}{3mm}
$
\begin{array}{|c|ccc|l|l|}
\hline
\text{No.} & n & a,b,c & \mu & \text{Link of $-\mu K_X-L$}
        & \text{other model}  \\
\hline
\hline
1 & 1 & 0,0,1 & 3 & \text{9-flop then $(2,0)$ to}
        & Y^\prime_{3,3}\subset\p^5(1^5,2)    \\
 & & & & \text{$\frac{1}{2}(1,1,1)$ singularity}
        & \text{general in its family}    \\
\hline
2 & 0 & 0,1,1 & 1 & \text{3-flop} & \text{$dP_3$ fibration, same}    \\
 & & & & & \text{numerology as $X$}      \\
\hline
3 & -1 & 1,1,1 & 1 & \text{flop} & \text{conic bundle over }\p^2    \\
 & & & & & \text{with $\deg\Delta = 7$}    \\
\hline
4 & -2 & 1,1,2 & 1 & \text{flop then $(2,1)$ to} &
         \text{$Y^\prime_4\subset \p^4(1^4,2)$}      \\
 & & & & \text{linear }\p^1\cong\ell\subset Y^\prime
        &      \\
\hline
5 & -2 & 1,2,2 & 1 & \text{Francia antiflip}
        & \text{$dP_2$ fibration with}  \\
 & & & & \text{then flop} & \text{$\frac{1}{2}(1,1,1)$ on 1 fibre}      \\
\hline
6 & -3 & 1,2,3 & 1 & \text{Francia antiflip then} &
         Y^\prime_6\subset \p^4(1^3,2,3)   \\
 & & & & \text{$(2,0)$ to $P\in Y^\prime$}
         & \text{$P$ a $c\mathrm{D}_4$ singularity}   \\
\hline
7 & -3 & 1,3,3 & 1 & \text{toric antiflip}
         & \text{$dP_1$ fibration with}        \\
 & & & & (1,1,-1,-3) & \text{$\frac{1}{3}(1,1,2)$ on 1 fibre}       \\
\hline
\hline
8a & -1 & 1,1,2 & 5
        & (1,1,-1,-1,-3),
        & Y^\prime\subset\p(1^4,2,3,4) \\
 & & &
        & \text{7-flop}, (2,0)
        & \text{general}, P=\frac{1}{4}(1,1,3)      \\
\hline
8b & -1 & 1,1,2 & 3
        & (1,1,-1,-1,-4),
        & \text{$Y^\prime\subset\p^5(1^4,2^2)$\,?}      \\
 & & &
        & \text{3-flop, $(2,0)$}
        &       \\
\hline
9 & -2 & 1,1,3 & 3
        & (1,1,-1,-1,-4),
        & \text{$Y^\prime\subset\p^5(1^2,2^2,3,5)$\,?}  \\
 & & &
        & \text{3-flop, $(2,0)$}
        &       \\
\hline
10 & -2 & 1,2,3 & 3
        & (1,1,-1,-2,-7),
        & \text{$Y^\prime\subset\p^5(1^2,2^3,3)$\,?}   \\
        & & & & (2,0) &      \\
\hline
\end{array}
$
\end{center}
\caption{Nonrigid $dP_3$ fibrations $X\in |3M+nL|\subset\F(0,a,b,c)$
\label{tab:2}}
\end{table}

\subsubsection{Overview of Table~\ref{tab:2}}

Each entry of Table~\ref{tab:2} represents a family of $dP_3$ fibrations
\[
X\in |3M+nL| \subset \F(0,a,b,c)
\] 
for which the 2-ray game on some member $X$ results in a
Sarkisov link to another model of $X$ as a Mfs.  The general member $X$
of families~1--7 is nonsingular.  In families~3--7, every $X$ (not just
the general member) admits the link as described.
Families~8--10 necessarily have a singularity on $\Gamma$, which
is described (in new coordinates) for general $X$ as follows.
\begin{center}
\renewcommand{\arraycolsep}{3mm}
$
\begin{array}{|l|c|c|c|c|}
\hline
\text{No.} & 8a & 8b & 9 & 10   \\
\hline
\text{Equation} & xy=zt & xy = z^3 + t^3 & xy = z^3 + t^3 & xy = z^3 + t^6 \\
\hline
\end{array}
$
\end{center}
In families~8--10, only special members $X$ admit the link as
described.

Our method for calculating the other model is to calculate the
graded ring of $-\mu K_X - L$.  In these four cases, though,
this does not present the other model well, since, we believe,
the contraction is to a non-Gorenstein singularity $P\in Y^\prime$.
We have not yet made the required calculations, but it seems
likely that the other model in families 8b, 9, 10 is better
presented as a complete intersection in weighted projective
spaces $\p^5(1^4,2^2)$, $\p^5(1^2,2^2,3,5)$ and $\p^5(1^2,2^3,3)$
respectively, the index of $P\in Y^\prime$ being 2, 5, 2 respectively.

We draw attention to family~2.  This links pair of $dP_3$ fibrations
$X\dasharrow X'$ that lie in the same family.  We guess that, in general,
they are not isomorphic and that they form a bi-rigid pair.  Indeed,
if they were isomorphic, this link would be an `untwisting' link
(in the sense of \cite{cpr}), and we would be inclined to guess
that $X$ is rigid.  We know little about this case.

The numbering of cases in the first column is arbitrary.  Apart from
this, the information for each entry is separated into three columns,
and we describe the contents of each of these in turn:
\begin{enumerate}
\item
This lists the integers $n$ and $a,b,c$ that determine the family
in question, as well as the $\mu>0$ for which
$-\mu K_X - L$ determines an edge of the mobile cone $\NM^1(X)$.
Note that, as Conjecture~\ref{con:dp-ast} predicts, condition~$(\ast)$
is not satisfied by these examples
\item
This describes the link.  We express the antiflips as $\cstar$ quotients
by listing the characters of a $\cstar$ action; we comment further below.
The word `flop' means the flop of a single rational curve.  We say
`$n$-flop' when, for general $X$, an analytic neighbourhood of the
flopping curve consists of the disjoint union of $n$ rational curve
neighbourhoods.
The notation $(2,p)$ indicates a divisorial contraction
to a point (when $p=0$) or a line (when $p=1$).
\item
This final column gives a sketch of the other model.
\end{enumerate}
We say more about the antiflips.  The Francia antiflip replaces a
$\p^1$ in the nonsingular locus having normal bundle
${\cal O}(-1)\oplus{\cal O}(-2)$ by a $\p^1$ passing through an
index~2 terminal quotient singularity.  This is worked out in
Section~\ref{sec:family5}.  The notation $(1,1,-1,-3)$ denotes a
toric \mbox{3-fold} antiflip similar to the Francia flip, but with an
index~3 singularity.  An example of the hypersurface antiflips
$(1,1,-1,-1,-m)$ is described in Section~\ref{sec:family8}:
typically, these replace a single rational curve passing through a
terminal Gorenstein point with a bouquet of $m-1$ rational curves
meeting in an index~$m$ singular point.  See \cite{MR2000f:14018}
for details and lists of these flips.

\paragraph{Other models as strict Mfs}
Four cases link to another strict Mfs.  Family~2 is discussed above,
and family~3 is linked to the conic bundle of Section~\ref{sec:cb7-2}.

The other two, families~5 and 7, share two novelties.
First, the use of weighted scrolls
for describing $dP_1$ and $dP_2$ fibrations contrasts with Grinenko's
use of finite morphisms to nonsingular scrolls.
As calculated in Section~\ref{sec:family5}, the other model of family~5
is the general element
\[
X^\prime\in |4M-L|
\quad
\text{in the weighted scroll}
\quad
\left(
\begin{array}{rrrrrr}
0 & 0 & 1 & 2 & 1 & 1   \\
1 & 1 & 0 & -1 & -1 & -1
\end{array}
\right).
\]
Similarly, one sees the other model of family~7 as the general element
\[
X^\prime\in |6M-3L|
\quad
\text{in the weighted scroll}
\quad
\left(
\begin{array}{rrrrrr}
0 & 0 & 2 & 3 & 1 & 1   \\
1 & 1 & -1 & -2 & -1 & -1
\end{array}
\right).
\]
Second, these examples present nonrigid $dP_1$ and $dP_2$ fibrations.
Grinenko \cite{gri02} has complete classifications of such nonrigid $dP_k$
in the Gorenstein case, and the nonrigid examples are rare.  The specimens
here have singularities of index~2 and 3, so they are not subject to
Grinenko's classifications, but they do invite one to extend Grinenko's
results to the higher index case.

\subsubsection{Family~5: general members are nonrigid}
\label{sec:family5}

We work out the link arising from the 2-ray game for
family~5 using the wall-crossing methods of Section~\ref{sec:cb7-3}
and the Appendix.  The calculations are similar for families~2--7.
Consider any $dP_3$ fibration
\[
X\colon (F=0) \in |3M - 2L| \subset \F = \F(0,1,2,2).
\]
The polynomial $F$ is a combination of the monomials of the Newton polygon
\[
\begin{array}{cc}
\text{$\deg_{u,v}$} & \text{fibre monomial}   \\
\hline
0 & xy^2, x^2z, x^2t    \\
1 & y^3, xyz, xyt       \\
2 & y^2z, y^2t, xz^2, xzt, xt^2 \\
3 & yz^2, yzt, yt^2     \\
4 & z^3, z^2t, zt^2, t^3
\end{array}
\]
where $\deg_{u,v}$ denotes the degree in $u,v$ of the coefficient.
We know that $F$ must involve both $xy^2+\varphi_1(u,v)y^3$
and $x^2z + x^2t$ nontrivially.
We continue to use the notation $\Gamma\colon (y=z=t=0)\subset X$.

The chamber defining $\F$ is $\r_+[L]+\r_+[M]$.  Crossing the wall
into the next chamber $\r_+[M]+\r_+[M-L]$ corresponds to a birational
map $\F\dasharrow\F_1$ that factors through the contraction of
$\Gamma\subset\F$ and the extraction of a surface
$E_1\colon (u=v=0)\subset\F_1$.  In fact, $E_1\cong\p(1,2,2)$ with
coordinates $y,z,t$.

The birational image $X_1\subset\F_1$ of $X$ is still $F=0$ so
$X_1\cap E_1$ is a $\p(1,2)$: when $u=v=0$ we can set $x=1$, as usual
in projective geometry, so the intersection is $y^2+z+t=0$ in $\p(1,2,2)$.
We see that $-K\Gamma = (M-L)\Gamma = -1$, so the map
$X\dasharrow X_1$ is an antiflip of $\Gamma$ and the link proceeds
with $X_1$ if and only if $X_1$ has terminal singularities.
We check this condition in coordinates.  

We calculate one patch of $X_1$ in detail, and leave the others to
the reader.  The unstable locus of the quotient defining $\F_1$ is
\[
(u=v=x=0) \cup (y=z=t=0)
\]
so the open set $xz\not= 0$ is a well-defined affine patch $U$ on $\F_1$.
There is a residual $\z/2$ action (the stabiliser of $z$-axis by the
$\mu$-action) defining the chart $\a^4\to U$ that acts on coordinates
$u,v,y,t$ of $\a^4$ by the character $(1,1,1,0)$.  The equation of $X_1$
in $U$ includes the monomial $t$, and the equation of the antiflipped
curve $\Gamma_1\subset X_1$ is $u=v=t=0$.  An analytic neighbourhood of
$X_1$ near the origin in $U$ is isomorphic to $t=0$ in the quotient
$\a^4/(\z/2)(1,1,1,0)$, and this is terminal.
The flip $X_1\to X$ is the Francia flip.

We cross to the next chamber giving $\F_1\dasharrow\F^\prime$
which factors as the contraction of $\p^2\subset\F_1$ (with
coordinates $u,v,x$) and the extraction of a $\p^1\subset\F^\prime$
(with coordinates $z,t$).
The birational image of $X$ is $X^\prime\subset\F^\prime$.
Since $F$ involves $xy^2 + uy^3$, the intersection of $X$ with
the exceptional $\p^2$ is a line $\Sigma_1$.  Clearly $-K\Sigma_1=0$
so the map $X_1\dasharrow X^\prime$ is a flop and, in particular, $X^\prime$ has
terminal singularities.

The variety $\F^\prime$ has a morphism $\F^\prime\to\p^1$ given by the ratio
$z,t$; fibres are $\p(1,1,2,1)$ with coordinates $u,v,x,y$.
To see this, make row operations on the character matrix
(basis changes in $\c^\times\times\c^\times$ that do not alter the quotients):
\[
\left(
\begin{array}{rrrrrr}
0 & 0 & 1 & 1 & 1 & 1     \\
1 & 1 & 0 & -1 & -2 & -2
\end{array}
\right)
\sim
\left(
\begin{array}{rrrrrr}
-1 & -1 & -1 & 0 & 1 & 1   \\
1 & 1 & 2 & 1 & 0 & 0
\end{array}
\right).
\]
The right-hand matrix is better adapted to the unstable locus
$(u=v=x=y=0) \cup (z=t=0)$ for $\F^\prime$.  The appropriate basis of
$\Pic(\F^\prime)$ is $M^\prime$, $L^\prime$, the line bundles corresponding
to the characters $0\choose 1$ and $1\choose 0$ in the basis of the
right-hand matrix.

The induced map $X^\prime\to\p^1$ describes $X^\prime$ as a $dP_2$ fibration:
since $x^2t\in F$, it is clear that
\[
X^\prime\subset |4M^\prime - L^\prime| \subset \F^\prime.
\]
One can check that $X^\prime$ is a general in its family if $X$ was
to start with.

\paragraph{}
In this example, every step of the 2-ray game on $X$ was inherited
(by computing birational images) from the steps of the 2-ray game on
$\F$, itself an easy toric calculation.  Most examples we know
have this feature.
\begin{dfn}
Let $X\subset\F$ be a $dP_3$ fibration in a scroll $\F\to\p^1$ (of
any dimension).
We say that \emph{the link follows the scroll} if the birational
images of $X$ in the 2-ray game of $\F$, together with the birational maps
induced between them, make up the 2-ray game of $X$.
\end{dfn}
This property is one of the characteristics of \emph{Mori dream spaces},
as introduced by Keel--Hu \cite{MR2001i:14059}.  Our less precise
definition is a convenient shorthand for the purposes of this paper only.

\subsubsection{Family~8: special members are nonrigid}
\label{sec:family8}

Consider $dP_3$ fibrations in the family
\[
X\colon (F=0) \in |3M-L|\subset\F = \F(0,1,1,2).
\]
The polynomial $F$ is a combination of the monomials of the Newton polygon
\[
\begin{array}{cc}
\text{$\deg_{u,v}$} & \text{fibre monomial}   \\
\hline
0 & x^2y, x^2z  \\
1 & xy^2, xyz, xz^2, x^2t       \\
2 & y^3, y^2z, yz^2, z^3, xyt, xzt      \\
3 & y^2t, yzt, z^2t, xt^2       \\
4 & yt^2, zt^2  \\
5 & t^3
\end{array}.
\]
We define
\[
{\mathrm{val}}(F) =
\min\{ \deg_{u,v} m\ |\ m\in F \text{ a monomial with nonzero coefficient} \}
\]
to separate cases of $F$ according to their leading term in $u,v$.

\paragraph{Case 1: ${\mathrm{val}}(F) = 0$.}

This is the general case: the coefficient of $x^2y$ or $x^2z$ in $F$ is
not zero and so the birational link of such $X$ follows the scroll.
Since $-K_X=M-L$, one calculates that the final divisorial contraction
is trivial against $-K_X$, so this is a bad link on $X$.

\paragraph{Case 2: ${\mathrm{val}}(F) = 1$.}

Every term of $F$ is divisible by $u$ or $v$ so we write
\[
F=uf - vg
\quad
\text{with}
\quad
f,g\in |3M-2L|.
\]
Since $u$ and $v$ do not vanish simultaneously on $X$, the rational
section
\begin{equation}
\label{eq:unproj}
\xi = f/v = g/u
\quad
\text{of}
\quad
{\cal O}_X(3M-3L)
\end{equation}
is regular: $\xi\in H^0(X,3M-3L)$.  We regard $\xi$ as a new variable
and, using equations (\ref{eq:unproj}), recompute the link starting from
\[
X \colon
(v\xi = f, u\xi = g)
\in |3M-2L| \cap |3M-2L| \subset
\F^5
\]
where the weighted scroll $\F^5$ is defined by the $\cstar\times\cstar$ action
\[
\left(
\begin{array}{rrrrrrr}
0 & 0 & 1 & 1 & 1 & 3 & 1       \\
1 & 1 & 0 & -1 & -1 & -3 & -2
\end{array}
\right)
\]
on $\a^7$ with coordinates $u,v,x,y,z,\xi,t$.  We have re-embedded $X$
isomorphically in a larger scroll: projection from $\xi$ is the
isomorphism to the original embedding.  Necessarily $F\owns x^2t$, since
otherwise $X$ would be singular along $\Gamma$ (the negative section
$y=z=\xi=t=0$), and then either $f\owns x^2t$, or $g\owns
x^2t$. Possibly after renaming $u,v$, we may assume that $f\owns x^2t$.
This case now divides into two subcases.

\paragraph{Case 2a: ${\mathrm{val}}(F) = 1$, ${\mathrm{val}}(g) = 0$.}

We can write the equations of $X$ in the form
\begin{equation}
\label{eq:2a}
\left(
\begin{array}{ccc}
 f_2 & f_1 + \xi & f_3  \\
 g_1 + \xi & g_2 & g_3  \\
\end{array}
\right)
\left(
\begin{array}{c}
 u \\ v \\ x
\end{array}
\right)
= 0.
\end{equation}
Following the link of $\F^5$, we see $X\dasharrow X_1$.  We solve
(\ref{eq:2a}) by 
\begin{equation}
\label{eq:eta}
\eta =
\frac{f_2g_2 - (f_1+\xi)(g_1+\xi)}{x} =
\frac{f_2g_3-f_3(g_1+\xi)}{v} =
\frac{(f_1+\xi)g_3-f_3g_2}{u}
\end{equation}
and conclude that $\eta \in H^0(X_1,5M-6L)$ since the semistable locus of
the action defining $X_1$ does not include $u=v=x=0$.  Once more, we
make a new scroll $\F^6$, by including $\eta$ among the coordinates,
and re-embed $X\subset\F^6$:
\[
X \colon ({\mathrm{Pfaff}}_{4\times 4}M = 0) \subset
\F^6
\]
where $\F^6$ is defined by the $\cstar\times\cstar$ action
\[
\left(
\begin{array}{rrrrrrrr}
0 & 0 & 1 & 1 & 1 & 3 & 5 & 1   \\
1 & 1 & 0 & -1 & -1 & -3 & -6 & -2
\end{array}
\right)
\]
on $\a^8$ with coordinates $u,v,x,y,z,\xi,\eta,t$.  Manipulating
(\ref{eq:eta}), the equations of $X\subset\F^6$ are the $4\times 4$
Pfaffians of a $5\times 5$ skew symmetric matrix
\[
M =
\left(
\begin{array}{cccc}
\eta &  f_2 & f_1 + \xi & f_3  \\
 & g_1 + \xi & g_2 & g_3  \\
 & & x & -v     \\
 & & & u
\end{array}
\right).
\]
The 2-ray game of $X$ follows the link of the scroll $\F^6$.
We sketch the steps:
\begin{enumerate}
\item
Antiflip $X\dasharrow X_1$.  This factors as $X\to Z_1\leftarrow X_1$,
by the contraction of $\Gamma$ to $Z_1$, both maps being the restrictions
of maps of the scroll.  The exceptional locus $u=v=0$ of
$Z_1\leftarrow X_1$ is defined by the equations of $X_1$ in
$\p^4(1,1,3,5,1)$ with coordinates $y,z,\xi,\eta,t$, since that $\p^4$ is
exceptional in the scroll.  We set $x=1$ to see that $t=0$, since
$f_3\owns xt$, and also that $\eta$ is eliminated.
So antiflipped locus is
\[
\bigl\{\text{three $\p(1,3)$s meeting in a point}\bigr\} =
\bigl(g_3(y,z) = 0\bigr)\subset\p(1,1,3).
\]
One can check in coordinates, as in Section~\ref{sec:family5}, that this
point is an index~3 terminal singularity on $X_1$.
In Table~\ref{tab:2}, this antiflip is denoted by $(1,1,-1,-1,-3)$.
\item
Flop $X_1\dasharrow X^\prime$.  Since $-K_{X^\prime}=M-L$, curves
that are contracted are trivial against $-K_{X^\prime}$.
See \cite{BZ03} for a detailed calculation, that also counts the
number of contracted curves.
\item
Divisorial contraction to a point $X^\prime\to Y^\prime$.
The divisor $F=(t=0)$ on the scroll is contracted to a point, so
$F\cap X^\prime$ is too.  The linear system $|6M-5L|$ (and its multiples)
define the morphism, so since
\[
-6K_X = 6(M-L) = (6M-5L) + L
\]
is negative on only the flipped curves, it must be relatively ample.
Therefore $X^\prime\to Y^\prime$ is extremal, and $Y^\prime$ has
terminal singularities.
\end{enumerate}
The result $Y^\prime$ is No.\ 6 of  Alt{\i}nok's list of codimension~3
Fano \mbox{3-folds} \cite{alt98}.  This Sarkisov link is already known to us,
calculated from the $Y^\prime$ end of the link in \cite{BZ03} following
Example~9.16 of \cite{kinosaki}.

\paragraph{Case 2b: ${\mathrm{val}}(F) = 1$, ${\mathrm{val}}(g) = 1$.}

Now $g$ is further specialised:
\[
g=uf_1 + vg_1
\quad
\text{with}
\quad
f_1,g_1\in |3M-3L|.
\]
As usual, since $u,v$ do not vanish simultaneously on $X$ and $u\xi - g$
is identically zero we write
\begin{equation}
\label{eq:unproj2}
u\xi = uf_1 + vg_1
\quad
\text{and}
\quad
\eta = (\xi - f_1)/v = g_1/u
\end{equation}
and conclude that $\eta\in H^0(X,3M-4L)$.  Again we include $\eta$ as
a new variable and compute the link again starting from
\begin{multline*}
X \in |3M-2L| \cap |3M-3L| \subset \F^5 \\
\text{determined by the action}
\quad
\left(
\begin{array}{rrrrrrr}
0 & 0 & 1 & 1 & 1 & 3 & 1       \\
1 & 1 & 0 & -1 & -1 & -4 & -2
\end{array}
\right).
\end{multline*}
Unlike case 2a, here we have eliminated $\xi$ from the coordinates
using the equation $\xi = v\eta + f_1$ of (\ref{eq:unproj2}).
Substituting for $\xi$ calculates the equations of $X\subset\F^5$:
\[
v^2\eta + vf_1 = f,\quad
u\eta = g_1.
\]
The 2-ray game of $X$ follows the link of the scroll $\F^5$.
The final contraction in the link is given by the linear system
$|3M-4L|$ and its multiples on $X^\prime$.

\subsubsection{Unstable $dP_3$ fibrations and condition $(\ast)$}
\label{sec:unstable}

Consider the family of $dP_3$ fibrations
\[
X\colon (F=0) \in |3M-4L|\subset\F(0,2,2,4).
\]
A general $X$ has a 2-ray game that follows the scroll, but this
gives a $K_X$-trivial bad link.  Nevertheless, we find
a special $X$ which is nonrigid.

The polynomial $F$
is a combination of the monomials of the Newton polygon
\[
\begin{array}{cc}
\text{degree} & \text{fibre monomial}   \\
\hline
0 & xy^2, xyz, xz^2, x^2t  \\
2 & y^3, y^2z, yz^2, z^3, xyt, xzt      \\
4 & y^2t, yzt, z^2t, xt^2       \\
6 & yt^2, zt^2  \\
8 & t^3 
\end{array}
\]
and we impose conditions on the coefficient polynomials $\alpha(u,v)$:
we require $u^i$ to divide the coefficient of fibre monomials according
to the table
\[
\begin{array}{cc}
i & \text{fibre monomial}   \\
\hline
1 & xyt, xzt      \\
2 & y^2t, yzt, z^2t     \\
3 & xt^2       \\
4 & yt^2, zt^2  \\
6 & t^3 
\end{array}
\]
with coefficients otherwise general.  The reader can check that a
general such $X$ is nonsingular away from a $c{\mathrm{D}}_4$
singularity at the point $(0,1;0,0,0,1)$.
Note that $X$ satisfies condition~$(\ast)$ because, even though it is
special in the family, the defining equation of $X$ involves $x^2t$ 
with nontrivial coefficient, hence it is still true that the 2-ray
game from $X$ follows the scroll (and ends in a bad link). 

Following Corti--Koll\'ar \cite{MR98e:14037}, \cite{MR98g:11076}, 
$X$ is unstable with respect to the
weight system $w = (3,2,2,0)$.  Indeed, $F$ is divisible by $u^6$
after the substitution
\[
u^3x^\prime, u^2y^\prime, u^2z^\prime, t^\prime
\quad
\text{for}
\quad
x,y,z,t.
\]
Cancelling the $u^6$ factor gives a $dP_3$ fibration
\[
X_{\mathrm{st}}\in |3M-L|\subset\F(0,1,1,1)
\]
that is square birational to $X\to\p^1$. Note that the fibre at $u=0$
has an Eckardt point:
\[
X_{\mathrm{st}}\cap (u=0)  = \bigl(tf_2(x,y,z) = g_3(y,z)\bigr).
\]
Even though it is not general in its family, $X_{\mathrm{st}}$ 
fails condition~$(\ast)$, and $X_{\mathrm{st}}$ does have a Sarkisov link,
following the scroll, to a conic bundle over $\p^2$.

\addcontentsline{toc}{section}{References}

\bibliography{bibbcz}

\def\cprime{$'$} \def\cprime{$'$} \def\cprime{$'$} \def\cprime{$'$}
  \def\cprime{$'$} \def\cprime{$'$} \def\cprime{$'$} \def\cprime{$'$}
\begin{thebibliography}{KMMT00}

\bibitem[ABR02]{MR1941620}
Selma Alt{\i}nok, Gavin Brown, and Miles Reid.
\newblock Fano 3-folds, {$K3$} surfaces and graded rings.
\newblock In {\em Topology and geometry: commemorating SISTAG}, volume 314 of
  {\em Contemp. Math.}, pages 25--53. Amer. Math. Soc., Providence, RI, 2002.

\bibitem[Alt98]{alt98}
S.~Alt{\i}nok.
\newblock {\em Graded rings corresponding to polarised {$K3$} surfaces and
  {$\mathbb{Q}$}-{F}ano 3-folds}.
\newblock PhD thesis, University of Warwick, 1998.

\bibitem[AM72]{MR48:299}
M.~Artin and D.~Mumford.
\newblock Some elementary examples of unirational varieties which are not
  rational.
\newblock {\em Proc. London Math. Soc. (3)}, 25:75--95, 1972.

\bibitem[Bro99]{MR2000f:14018}
Gavin Brown.
\newblock Flips arising as quotients of hypersurfaces.
\newblock {\em Math. Proc. Cambridge Philos. Soc.}, 127(1):13--31, 1999.

\bibitem[BZ]{BZ03}
Gavin Brown and Francesco Zucconi.
\newblock The graded ring of a rank 2 {S}arkisov link.
\newblock In preparation.

\bibitem[Cat81]{MR83c:14026}
F.~Catanese.
\newblock Babbage's conjecture, contact of surfaces, symmetric determinantal
  varieties and applications.
\newblock {\em Invent. Math.}, 63(3):433--465, 1981.

\bibitem[CM]{cm02}
Alessio Corti and Massimiliano Mella.
\newblock Birational geometry of terminal quartic \mbox{3-folds}.\ {I}.
\newblock to appear in Amer.\ Jour.\ Math.

\bibitem[Cor95]{MR96c:14013}
Alessio Corti.
\newblock Factoring birational maps of threefolds after {S}arkisov.
\newblock {\em J. Algebraic Geom.}, 4(2):223--254, 1995.

\bibitem[Cor96]{MR98e:14037}
Alessio Corti.
\newblock Del {P}ezzo surfaces over {D}edekind schemes.
\newblock {\em Ann. of Math. (2)}, 144(3):641--683, 1996.

\bibitem[Cor00]{corti:00}
Alessio Corti.
\newblock Singularities of linear systems and {$3$}-fold birational geometry.
\newblock In {\em Explicit birational geometry of 3-folds}, volume 281 of {\em
  London Math. Soc. Lecture Note Ser.}, pages 259--312. Cambridge Univ. Press,
  Cambridge, 2000.

\bibitem[CPR00]{cpr}
Alessio Corti, Aleksandr Pukhlikov, and Miles Reid.
\newblock Fano {$3$}-fold hypersurfaces.
\newblock In {\em Explicit birational geometry of 3-folds}, volume 281 of {\em
  London Math. Soc. Lecture Note Ser.}, pages 175--258. Cambridge Univ. Press,
  Cambridge, 2000.

\bibitem[CR00]{MR2001f:14004}
Alessio Corti and Miles Reid, editors.
\newblock {\em Explicit birational geometry of \mbox{3-folds}}, volume 281 of
  {\em London Mathematical Society Lecture Note Series}.
\newblock Cambridge University Press, Cambridge, 2000.

\bibitem[Dix02]{Dixon}
A.~C. Dixon.
\newblock Note on the reduction of a ternary quantic to a symmetric
  determinantal.
\newblock {\em Proc. Cam. Phil. Soc.}, 1902.

\bibitem[Gria]{gri02}
M.~M. Grinenko.
\newblock Birational models of del {P}ezzo fibrations.
\newblock math.AG/0209394.

\bibitem[Grib]{grin03}
M.~M. Grinenko.
\newblock Non-rationality of a three-dimensional {F}ano variety of index 2 and
  degree 1.
\newblock math.AG/0304321.

\bibitem[Gric]{grin99}
M.~M. Grinenko.
\newblock On a rigidity criterion for del {P}ezzo fibrations over
  {$\mathbb{P}^1$}.
\newblock math.AG/9911210.

\bibitem[Gri00a]{grinenko:00}
M.~M. Grinenko.
\newblock Birational properties of pencils of del {P}ezzo surfaces of degrees 1
  and 2.
\newblock {\em Mat. Sb.}, 191(5):17--38, 2000.

\bibitem[Gri00b]{grinenko:01c}
M.~M. Grinenko.
\newblock On the birational rigidity of some pencils of del {P}ezzo surfaces.
\newblock {\em J. Math. Sci. (New York)}, 102(2):3933--3937, 2000.
\newblock Algebraic geometry, 10.

\bibitem[Gri01a]{grinenko:01a}
M.~M. Grinenko.
\newblock On fiberwise surgeries of fibrations on a del {P}ezzo surface of
  degree 2.
\newblock {\em Uspekhi Mat. Nauk}, 56(4(340)):145--146, 2001.

\bibitem[Gri01b]{grinenko:01b}
M.~M. Grinenko.
\newblock On fibrations into del {P}ezzo surfaces.
\newblock {\em Mat. Zametki}, 69(4):550--565, 2001.

\bibitem[HK00]{MR2001i:14059}
Yi~Hu and Sean Keel.
\newblock Mori dream spaces and {GIT}.
\newblock {\em Michigan Math. J.}, 48:331--348, 2000.
\newblock Dedicated to William Fulton on the occasion of his 60th birthday.

\bibitem[IP96]{MR98e:14009}
V.~A. Iskovskikh and A.~V. Pukhlikov.
\newblock Birational automorphisms of multidimensional algebraic manifolds.
\newblock {\em J. Math. Sci.}, 82(4):3528--3613, 1996.
\newblock Algebraic geometry, 1.

\bibitem[Isk87]{MR88i:14038}
V.~A. Iskovskikh.
\newblock On the rationality problem for conic bundles.
\newblock {\em Duke Math. J.}, 54(2):271--294, 1987.

\bibitem[Isk91a]{MR93i:14033}
V.~A. Iskovskih.
\newblock Towards the problem of rationality of conic bundles.
\newblock In {\em Algebraic geometry (Chicago, IL, 1989)}, volume 1479 of {\em
  Lecture Notes in Math.}, pages 50--56. Springer, Berlin, 1991.

\bibitem[Isk91b]{MR92k:14036}
V.~A. Iskovskikh.
\newblock On the rationality problem for conic bundles.
\newblock {\em Mat. Sb.}, 182(1):114--121, 1991.

\bibitem[Isk95]{MR2000i:14019}
V.~A. Iskovskikh.
\newblock On the rationality problem for three-dimensional algebraic varieties
  fibered over del {P}ezzo surfaces.
\newblock {\em Trudy Mat. Inst. Steklov.}, 208(Teor. Chisel, Algebra i Algebr.
  Geom.):128--138, 1995.
\newblock Dedicated to Academician Igor\cprime\ Rostislavovich Shafarevich on
  the occasion of his seventieth birthday (Russian).

\bibitem[Isk96]{MR97j:14044}
V.~A. Iskovskikh.
\newblock On a rationality criterion for conic bundles.
\newblock {\em Mat. Sb.}, 187(7):75--92, 1996.

\bibitem[KMMT00]{MR2001h:14053}
J{\'a}nos Koll{\'a}r, Yoichi Miyaoka, Shigefumi Mori, and Hiromichi Takagi.
\newblock Boundedness of canonical {$\mathbb Q$}-{F}ano 3-folds.
\newblock {\em Proc. Japan Acad. Ser. A Math. Sci.}, 76(5):73--77, 2000.

\bibitem[Kol97]{MR98g:11076}
J{\'a}nos Koll{\'a}r.
\newblock Polynomials with integral coefficients, equivalent to a given
  polynomial.
\newblock {\em Electron. Res. Announc. Amer. Math. Soc.}, 3:17--27
  (electronic), 1997.

\bibitem[KSC03]{KSC}
J.~Koll\'{a}r, K.~Smith, and A.~Corti.
\newblock {\em Rational and nearly rational varieties}.
\newblock Cambridge University Press, 2003.

\bibitem[Mel]{mella02}
Massimiliano Mella.
\newblock Birational geometry of terminal quartic \mbox{3-folds}.\ {II}.
\newblock {S}ubmitted.

\bibitem[Puk97]{MR98j:14014}
A.~V. Pukhlikov.
\newblock Birational automorphisms of three-dimensional algebraic manifolds
  stratified on cubic surfaces.
\newblock {\em Uspekhi Mat. Nauk}, 52(1(313)):235--236, 1997.

\bibitem[Puk98a]{MR2000d:14017}
A.~V. Pukhlikov.
\newblock Birational automorphisms of three-dimensional algebraic varieties
  fibered by {D}el {P}ezzo surfaces.
\newblock {\em Dokl. Akad. Nauk}, 361(1):14--16, 1998.

\bibitem[Puk98b]{MR99f:14016}
A.~V. Pukhlikov.
\newblock Birational automorphisms of three-dimensional algebraic varieties
  with a pencil of del {P}ezzo surfaces.
\newblock {\em Izv. Ross. Akad. Nauk Ser. Mat.}, 62(1):123--164, 1998.

\bibitem[Puk00]{MR2001j:14010}
Aleksandr~V. Pukhlikov.
\newblock Essentials of the method of maximal singularities.
\newblock In {\em Explicit birational geometry of 3-folds}, volume 281 of {\em
  London Math. Soc. Lecture Note Ser.}, pages 73--100. Cambridge Univ. Press,
  Cambridge, 2000.

\bibitem[Rei]{kinosaki}
Miles Reid.
\newblock Graded rings and birational geometry.
\newblock In {\em Proc.\ of algebraic geometry symposium (Kinosaki, Oct.\
  2000)}.

\bibitem[Rei97]{MR98d:14049}
Miles Reid.
\newblock Chapters on algebraic surfaces.
\newblock In {\em Complex algebraic geometry (Park City, UT, 1993)}, volume~3
  of {\em IAS/Park City Math. Ser.}, pages 3--159. Amer. Math. Soc.,
  Providence, RI, 1997.

\bibitem[Sar79]{MR81d:14023}
V.~G. Sarkisov.
\newblock Birational automorphisms of three-dimensional algebraic varieties
  representable as a conic bundle.
\newblock {\em Uspekhi Mat. Nauk}, 34(4(208)):207--208, 1979.

\bibitem[Sar80]{MR82g:14035}
V.~G. Sarkisov.
\newblock Birational automorphisms of conic bundles.
\newblock {\em Izv. Akad. Nauk SSSR Ser. Mat.}, 44(4):918--945, 974, 1980.

\bibitem[Sar82]{MR84h:14047}
V.~G. Sarkisov.
\newblock On conic bundle structures.
\newblock {\em Izv. Akad. Nauk SSSR Ser. Mat.}, 46(2):371--408, 432, 1982.

\bibitem[Sob02]{MR2003g:14017}
I.~V. Sobolev.
\newblock Birational automorphisms of a class of varieties fibered by cubic
  surfaces.
\newblock {\em Izv. Ross. Akad. Nauk Ser. Mat.}, 66(1):203--224, 2002.

\bibitem[Tak02]{MR1924722}
Hiromichi Takagi.
\newblock On classification of {$\mathbb Q$}-{F}ano 3-folds of {G}orenstein
  index 2. {I}, {II}.
\newblock {\em Nagoya Math. J.}, 167:117--155, 157--216, 2002.

\end{thebibliography}

\appendix

\section{Birational transformations of scrolls}
\label{sec:scrolls}

\subsection{Definition of scrolls}
\label{sec:def-scrolls}

We set our notation for rational scrolls and toric links between them. Our
treatment follows closely \cite[Chapter 2]{MR98d:14049}.

Throughout the Appendix, we consider actions of the group
$G=\c^\times \times \c^\times$ on affine space. The elements of
$G$ are ordered pairs $(\lambda, \mu)$ where $\lambda, \mu \in
\c^\times$. We denote by $\mathbb{X}=\Hom (G, \c^\times)$
the lattice of characters of $G$, with basis the coordinate functions
(projections on the two factors) $\chi_1$, $\chi_2$ such that
$\chi_1(\lambda, \mu) = \lambda$ and $\chi_2 (\lambda, \mu) = \mu$. 
The dual lattice $\mathbb{X}^\ast = \Hom (\mathbb{X}, \z)$ is based by
the 1-parameter subgroups $e_1$, $e_2$ such that $e_1(\lambda) =
(\lambda, 1)$ and $e_2(\mu) = (1, \mu)$. Sometimes we abuse notation
and write $\lambda$ for the element $e_1(\lambda) \in
G$ (and, similarly, $\mu$ for $e_2(\mu)$). Occasionally we abuse even
further and identify $\lambda$ with the coordinate function
$\chi_1\colon G \to \c^\times$ (and, similarly, we identify 
$\mu$ with $\chi_2$). 

We now define rational scrolls.
Fix a base $\p=\p^k$, with homogeneous coordinates $u_0, \dots,
u_k$; in this paper, we only work with $k=1$ or $k=2$. Consider 
now $\c^{n+1}$, with 
coordinates $x_0,\dots,x_n$. Fix integers $a_0,\dots,a_n$, usually 
nonnegative and in increasing order. Consider the action of $G$ on the 
affine space $\a = \c^{k+1} \times \c^{n+1}$, where the two factors
of $G$ act by
\begin{align*}
 \lambda \colon (u_0,\dots,u_k,x_0,\dots,x_n) \mapsto & (\lambda u_0,\dots,
 \lambda u_k,\lambda^{-a_0} x_0,\dots,\lambda^{-a_n} x_n)\\
 \mu \colon (u_0,\dots,u_k,x_0,\dots,x_n) \mapsto & (u_0,\dots,u_k,\mu
 x_0,\dots,\mu x_n).
\end{align*}
We summarise this action by writing down the matrix:
\[
\begin{pmatrix}
  1 &\dots& 1 &-a_0&\dots&-a_n\\
  0 &\dots& 0 &1   &\dots&1
\end{pmatrix}.
\] 
By definition, the scroll $\F=\F(a_0,\dots,a_n)$ is the following quotient:
\[
\F(a_0,a_1,\dots,a_n) = \bigl(\c^{k+1}\setminus \{0\}\bigr) \times
\bigl(\c^{n+1}\setminus \{0\}\bigr)/G.
\]
It is clear that $\F$ is a $\p^n$-bundle over $\p^k$.

\subsection{Line bundles on scrolls}
\label{sec:lbs-scrolls}

There is a 1-to-1 correspondence between line bundles on the scroll
$\F$ and characters $\chi\colon G \to \cstar$ of the group
$G$. To a character $\chi$ we associate a
\emph{$G$-linearisation} of the trivial line bundle over $\a$, by
acting with $\chi$ in the direction of fibres. Taking the quotient of
the $G$-linearisation by $G$, we form a bundle on $\F$.
Let us denote by $L_\chi$ the resulting bundle. It is easy to
chase through the definition and see that the 
sections of $L_\chi$ are the global  
eigenfunctions with eigenvalue $\chi$:
\[
H^0(\F, L_\chi)=\{f\colon (\c^{k+1} \setminus \{0\}) \times 
(\c^{n+1} \setminus \{0\}) \to \c \mid f(gx)=\chi(g) f(x) \}. 
\]
Sometimes, we abuse notation and confuse the line bundle and the
corresponding character. 

We denote by $L$ and $M$ the line bundles corresponding to $\chi_1$ and
$\chi_2$. By what we just said, the 
sections of $L$ are the functions $f\colon (\c^{k+1}\setminus\{0\}) 
\times (\c^{n+1}\setminus \{0\}) \to \c$ such that
\begin{eqnarray*}
   f(\lambda u_0,\dots, \lambda u_k,\lambda^{-a_0} x_0,\dots,\lambda^{-a_n}
  x_n)&=& \lambda f(u_0,\dots,u_k,x_0,\dots,x_n)\\
  f(u_0,\dots,u_k,\mu x_0,\dots,\mu x_n)&=& f(u_0,\dots,u_k,x_0,\dots,x_n).
\end{eqnarray*}
It follows that $H^0(\F, L)$ is based by the coordinate
functions $u_0,\dots,u_n$ and $L$ is the pull-back $\pi^\ast \o(1)$ by
the natural morphism $\pi \colon \F\to \p^k$. 
Similarly, the 
group $H^0(\F, M)$ of global sections of $M$ is based by the
monomials 
\[
u_0^{w_0}\cdots u_k^{w_k} x_i \quad \text{where} \quad w_0+\cdots +w_k = a_k.
\]

\subsection{$\F$ as a geometric quotient}

Using the language of $G$-linearisations, we can view $\F$ as a
\emph{geometric} quotient $\a^{k+n+2} /\!\!/ G$ in the sense of
Geometric Invariant Theory. If $G$ acts on $\a=\a^{k+n+2}$ as before and
$\chi\colon G \to \cstar$ is a character, then $G$ acts on functions
$f\in \o_\a$ by $gf(x)=\chi(g)f(g^{-1}x)$, and we denote $\o_\a^\chi$ the
invariants. The set of semistable points is by definition
\[
\a^{\text{ss}}_\chi=\{x\in \a \mid \exists f\in \o_\a^\chi, \;f(x) \not = 0\}. 
\]
The group $G$ acts with finite stabilisers (in fact, freely) on the
set of semistable points, and the geometric quotient is by definition
\[
\a /\!\!/G = \a^{\text{ss}}_\chi/G.
\] 

\subsection{Linearisations and geometric quotients}
\label{sec:geom-quot}

Different linearisations lead to different quotients. We state what is
going on and leave the elementary proofs to the reader. 
We say that a linearisation is \emph{useful} if the set of semistable points is
nonempty. To fix ideas, let us assume that $0=a_0\leq a_1 \leq \dots
\leq a_{n}$.
The cone of useful linearisations is the cone
\[
(\r_+[L]+\r_+[M-a_{n-1} L])\cap \mathbb{X}
\]
This cone is partitioned into chambers corresponding to different
geometric quotients. For example if $\chi \in \r_+[L] + \r_+[M]$, then 
$\a^{\text{ss}}_\chi=(\c^{k+1}\setminus \{0\}) \times (\c^{n+1} \setminus
\{0\})$ and we get our $\F$ back. However there are the other chambers
\[
\sigma_i = (\r_+[M-a_{i-1}L] + \r_+[M-a_{i}L]) \cap \mathbb{X},
\]
whenever $a_{i-1}<a_{i}$, and if $\chi \in \sigma_i$
\[
\a^{\text{ss}}_\chi = (\c^{k+1+i}\setminus \{0\})\times (\c^{n+1-i}
\setminus \{0\}).
\]
The corresponding quotient 
\[
\F_i= \a^{\text{ss}}_\chi /G
\]
for $\chi \in \sigma_i$ is birational to $\F$. The sequence of
birational maps $\F=\F_0 \dasharrow \F_1 \dasharrow \cdots$ is a 2-ray game
in the sense of \cite[Section 2.2]{corti:00}. For example if $a_1>0$, the move 
$\F\dasharrow \F_1$ is the antiflip of the section $\Gamma=
\{x_1=x_2=\dots =x_n=0\}$ which generates one of the two extremal
rays of $\NE (\F)$.

\subsection{Generalisations}
\label{sec:generalisations}

The above can be generalised slightly to the action of $G$ on
$\a^{n+1}$ given by the matrix
\[
\begin{pmatrix}
  \alpha_0 & \dots & \alpha_n \\
  \beta_0  & \dots & \beta_n
\end{pmatrix}.
\]
We use this notation in the text as a shorthand for the action
\begin{align*}
  x_0, \dots, x_n \mapsto \lambda^{\alpha_0}x_0, \dots,
  \lambda^{\alpha_n} x_n \\
 x_0, \dots, x_n \mapsto \mu^{\beta_0}x_0, \dots,
  \mu^{\beta_n} x_n.
\end{align*}
Here we assume:
\begin{enumerate}
\item the rows are linearly independent (so we get a faithful action
  of $G$) and all columns are nonzero,
\item all $\beta_i\geq 0$, and $\beta_n >0$,
\item the ratios $\alpha_i/\beta_i$ are in decreasing order.  
\end{enumerate}
As before we denote $L$, $M$ the $G$-linearisations corresponding to
the characters $\chi_1$, $\chi_2$. The cone of useful linearisations
\[
(\r_+[\alpha_0M-\beta_0L]+\r_+[\alpha_{n-1}M-\beta_{n-1}L]) \cap \mathbb{X}
\]
is partitioned into chambers
\[
\sigma_i =
(\r_+[\alpha_{i-1}M-\beta_{i-1}L]+\r_+[\alpha_{i}M-\beta_{i}L]).
\]
(when $\alpha_{i-1}/\beta_{i-1}> \alpha_i/\beta_i $).
Choosing a character $\chi \in \sigma_i$ gives a semistable locus
$\a^{\text{ss}}_\chi= (\c^{i+1}\setminus \{0\}) \times (\c^{n-i}\setminus
\{0\})$ and a geometric quotient $\F_i=\a^{\text{ss}}_\chi/G$.

\subsection{Alternative viewpoints}
\label{sec:altern-viewp}

There are at least two other points of view on birational maps of
scrolls.

We can identify the cone of useful $G$-linearisations with the mobile
cone $\NM^1(\F)$ of the scroll. From this point of view, the main focus
is the bigraded ring 
\[
\bigoplus_{e,n} H^0(\F, eM+nL).
\]

We can also view scrolls, and their generalisations, as special cases
of rank~2 toric varieties.

\end{document}